\documentclass[11pt]{article}

\usepackage{latexsym}
\usepackage[latin1]{inputenc}
\usepackage{amsfonts,amsmath,amssymb}
\usepackage{todonotes}

\usepackage{a4wide}

\author{Andrea Surroca Ortiz\footnote{The first version of this paper was partially supported by the Marie Curie IEF 025499  fellowship of the European Community and the second one by the Ambizione fund PZ00P2\_121962 of the Swiss National Science Foundation. This last version is partially supported by the EPSRC EP/N007956/1 grant at the University of Manchester.}
}

\title{\Large Conjectural estimates on the Mordell-Weil and Tate-Shafarevich
groups of an abelian variety}
\date{\today}

\newtheorem{thm}{\textbf{Theorem}}[section]
\newtheorem{conj}[thm]{\textbf{Conjecture}}
\newtheorem{lemma}[thm]{\textbf{Lemma}}
\newtheorem{remark}[thm]{\textbf{Remark}}
\newtheorem{prop}[thm]{\textbf{Proposition}}
\newtheorem{hyp}[thm]{\textbf{Hypothesis}}

\newcommand\rat{\mathbb{Q}}
\newcommand\C{\mathbb{C}}
\newcommand\Z{\mathbb{Z}}
\newcommand\enteros{\mathbb{Z}}
\newcommand\R{\mathbb{R}}

\newcommand\Kbarre{\overline{K}}
  
\newcommand\tors{\mathrm{tors}}
\newcommand\card{\mathrm{card}}

\newcommand\ord{\mathrm{ord}}
\newcommand\rk{\mathrm{rk}}
\newcommand\codim{\mathrm{codim}}

 \newcommand\Pic{\mathrm{Pic}}

\newcommand\Reg{\mathrm{Reg}}
\newcommand\Vol{\mathrm{Vol}}
\newcommand\Covol{\mathrm{Covol}}
\newcommand\Res{\mathrm{Res}}
\newcommand\Gal{\mathrm{Gal}}

\newcommand\vp{\mathfrak{p}}
\newcommand\vq{\mathfrak{q}}
\newcommand\va{\mathfrak{a}}

\newcommand\A{\mathcal{A}}
\newcommand\F{\mathcal{F}}
\newcommand\G{\mathcal{G}}

\newcommand\pp{\mathcal{L}}

\newcommand\ok{\mathcal{O}_{K}}
\newcommand\Sok{\mathrm{Spec}(\ok)}

\newcommand\DiffAK{\Omega_{A/K}^{1}}
\newcommand\DiffA{\Omega_{\mathcal{A}/\ok}^{g}}

\newcommand\eDiffA{\omega_{\mathcal{A}/\ok}}

 \newcommand\Ad{\check{A}}
 
\newcommand\espproj{\mathbf{P}}

\newcommand\esp{\hspace{0,2cm}}

\input cyracc.def \font\tencyr=wncyr10 \def\russe{\tencyr\cyracc}
\def\Sha{\text{\russe{Sh}}}

\begin{document}

\bibliographystyle{alpha}

\maketitle

{\begin{flushright}
{\it To my daughter} 
\end{flushright}}

\begin{quote}
\textbf{Abstract.}

{\small We consider an abelian variety defined over a number field. We give conditional bounds for the order of its Tate-Shafarevich group, as well as conditional bounds for the N\'eron-Tate height of generators of its Mordell-Weil group. The bounds are implied by strong but, henceforth, classical conjectures, such as the Birch and Swinnerton-Dyer conjecture and the functional equation of the $L$-series. In particular, we improve and generalise a result by D. Goldfeld and L. Szpiro on the order of the Tate-Shafarevich group, and extend a conjecture of S. Lang on the canonical height of a system of generators of the torsion-free part of the Mordell-Weil group. The method is an extension of the algorithm proposed by Yu. Manin for finding  a basis for the non-torsion rational points of an elliptic curve defined over the rationals.} 
\end{quote}

{\small 2000 Mathematics Subject Classification: 11G40, 11G50, 11G10, 11G05, 14G05, 11H50, 14G40.}

{\small Keys words: Abelian variety, $L$-functions, heights, Birch--Swinnerton-Dyer conjecture, Mordell-Weil group, Tate-Shafarevich group.}

\section{Introduction}\label{intro}

 The Mordell-Weil theorem states that the group of rational points on an abelian variety $A/K$ defined over a number field is finitely generated, and can thus be written as
 $$A(K) \simeq A(K)_{\tors} \oplus \Z P_1 \oplus \ldots \oplus \Z P_r,$$
 where $r = \rk (A(K))$ is called its rank, and $A(K)_{\tors}$ is the finite group of its torsion points. 
 
  While there exist results  on the torsion part, the torsion-free part remains less tractable.  Even in the particular case of an elliptic curve, there is no way, in general, to compute the rank or a set of generators of this group. 
  
  The proof of the Mordell-Weil theorem involves the Tate-Shafarevich group $\Sha (A/K)$ of $A/K$, 
  which measures the obstruction to the Hasse principle. In fact, a non-trivial element of $\Sha (A/K)$ corresponds to a homogenous space, which has $K_{v}$-rational points for every place $v$, but no $K$-rational points. Even if it is not easy to construct such a variety, it is still unknown, in the general case,  if $\Sha (A/K)$ is finite.

 For some applications, it would be sufficient to bound the {\it size} (e.g. cardinality, height, volume) of the invariants related to the variety.  
  
 In this article, we explore how  the canonical height of a well-chosen system of generators (which provides arithmetic information) could be bounded, as well as the order of the Tate-Shafarevic group of $A(K)$. 
 The bounds given here are not conjectured, but implied, by strong but classical conjectures. 
 We follow the approach of Manin, who proposed a conditional algorithm for finding a basis for the non-torsion rational points of an elliptic curve over $\rat$. 
 The method is based on the hypothesis that  the $L$-series of the elliptic curve satisfies a functional equation and the celebrated conjecture of Birch and Swinnerton-Dyer \cite{birch-swinnerton-dyer} (BSD-conjecture for short), which translates analytic information into algebraic and arithmetic information.    

We extend Manin's method to an abelian variety $A$ of arbitrary dimension, defined over an arbitrary number field $K$. 
 The bounds are given in terms of more tractable objects associated to the variety  and the number field. Precisely, our bounds depends on the  Faltings' height $h = h_{Falt}(A/K)$ of $A/K$ (which measures the arithmetic complexity of the variety),   
the absolute value $\F = |N_{K/\rat}\F_{A/K}|$ of the norm of the conductor (which gives information about the places of bad reduction),  the dimension $g$ of $A$, the Mordell-Weil rank $r = \rk(A(K))$, the degree $d = [K:\rat]$, and the absolute value $D_K$ of the discriminant of $K$.  
Moreover, the dependence of the bounds is explicit in all the parameters. 

Suppose that $A$ carries a polarisation $\phi_{\pp} : A \to \Ad$. The associated Néron-Tate height on $A(\Kbarre)$, 
 $$\hat{h}_{\pp}  = \hat{h}_{A,\pp} : A(\Kbarre) \to \R,$$ 
extends to a positive definite quadratic form on  $A(K) \otimes_{\Z} \R$. The associated pairing $<,>_{\pp}$ endows $A(K) \otimes_{\Z} \R \simeq \R^{r}$   with a structure of a euclidean space, and we can view $A(K)/A(K)_{\tors}$ as a lattice sitting inside this  $\R^{r}$. The regulator $\Reg_{\pp}(A/K): = \det (<P_i, P_j>_{\pp} )_{1 \leq i, j \leq r} \geq 0$  of $A/K$ relative to $\pp$, where $\{ P_i\}_{1 \leq i \leq r}$ is a basis for $A(K)/A(K)_{\tors}$, is the square of the volume of the fundamental domain for the lattice. Furthermore, we can  define a {\it canonical} regulator $\Reg (A/K)$, independent of $\pp$.

 From Manin's algorithm one could deduce a bound for the product of the canonical regulator and the Tate-Shafarevich group of $A(K)$. 
  On this topic, a beautiful analogy with the classical Brauer-Siegel's formula (relating the discriminant, the regulator and the class number of a number field) is developed by Hindry in \cite{hindry.mordell-weil}.
He formulates the following  conjecture: 
{\it For all $\epsilon >0,$  
 \begin{equation}\label{hindry'sconj}
  |\Sha(A/K)| \cdot \Reg(A/K) \ll H_{Falt}(A/K)^{1+\epsilon},
\end{equation}
 where $H_{Falt}(A/K) = e^{h_{Falt}(A/K)}$ and the implicit constants in the $\ll$ symbol depend on $K, g, \epsilon$ and $\rk(A(K))$.}\footnote{See \cite{marc-amilcar-brauer-siegel} for an unconditional  function field analogue, and \cite{griffon-BS-Legendre-EC} for an upper and a lower similar bounds for a family of elliptic curves over ${\bf F}_q(t)$.} 
 This extends a conjecture of Lang \cite[p. 99]{Lang-NTIII-Enciclopeadia} for elliptic curves over the field of the rational numbers: {\it Let $E$ be an elliptic curve defined over $\rat$, with minimal equation over $\enteros$ given by $y^2 = x^3 -\gamma_2 x - \gamma_3$, and $H(E) = \max\{ |\gamma_2|^3, |\gamma_3|^2\}$. Then  
  \begin{equation}\label{lang'sconj-Sha.Reg}
   |\Sha(E/\rat)| \cdot \Reg(E/\rat) \leq b_1 H(E)^{1/12} \cdot \F ^{\epsilon(\F)}\cdot b_2^r \cdot (\log \F)^r,
   \end{equation}
for $b_1$ and $b_2$ absolute constants and $\epsilon(\F)$ tends to $0$, when $\F$ tends to $\infty$.}

 Lang  modified the heuristic approach of Manin and also proposed the following conjecture \cite[Conjecture 3]{lang.conj.dio}:
{\it Let $E$ be an elliptic curve defined over $\rat$. We can find a basis $\{P_{1}, \ldots, P_{r}\}$ for the torsion-free part of $E(\rat)$ satisfying}
 \begin{equation}\label{lang'sconj}
 \max_{1\leq i \leq r}\hat{h}(P_i) \ll c^{\rk(E(\rat))^{2}}\cdot \F_{E/\rat}^{\epsilon(\F_{E/\rat})} \cdot (\log \F_{E/\rat})^{\rk(E(\rat))} \cdot e^{h_{Falt}(E/\rat)},
 \end{equation}
{\it  where $c$ is an absolute constant and  $\epsilon$ is a function which does not depend on the rank, and   $\epsilon(\F)$  tends to 0 as $\F$ tends to infinity.}

Furthermore, 
from the proof of the Weak Mordell-Weil theorem, we know that, for all $n \geq 1$, the $n$-torsion part of the Tate-Shafarevich group is finite. It is conjectured that the whole Tate-Shafarevich group is finite; the conjecture is known for certain elliptic curves with complex multiplication (\cite{rubin1987}) and certain modular elliptic curves (\cite{kolyvagin1989}). Goldfeld and Szpiro \cite{goldfeld-szpiro} suggested the following bound for the order of the Tate-Shafarevich group $\Sha(E/K)$ of an elliptic curve, in terms of the conductor:  
{\it Let $E$ be an elliptic curve defined over a field $K$, which can be a number field or a function field. Then, for every $\epsilon >0$,  
\begin{equation}\label{conj.g-sz.sha<cond}
|\Sha(E/K)| = O(\F_{E/K}^{1/2 + \epsilon}),
\end{equation}
where the implicit constant in the $O$ depends on $\epsilon$, $K$ and $\rk(E(K))$.}
In the same article, they announced that this conjecture 
holds for elliptic curves defined over function fields provided the Tate-Shafarevich group of the function field is finite. Independently, Rajan \cite{rajan} proved this result. 
Goldfeld and Lieman \cite{goldfeld-lieman} proved that for a CM elliptic curve defined over $\rat$ with Mordell-Weil rank 0, we have $|\Sha(E/\rat)| < k(\epsilon) \F_{E/\rat}^{\delta + \epsilon}$, with $\delta = \frac{59}{120}$ if $j \ne 0, 1728$, $\delta = \frac{37}{60}$ if $j = 0$ and $\delta = \frac{79}{120}$ if $j = 1728$, where $k(\epsilon)$ depends only on $\epsilon$ and is effectively computable.
 It is also proved in \cite[Theorem 1]{goldfeld-szpiro} that, {\it if the curve $E$ is defined over $\rat$ and satisfies the BSD-conjecture and Szpiro's conjecture} (which  predicts a bound for the discriminant in terms of the conductor), {\it then} 
 \begin{equation}\label{thm1-Goldfeld-Szpiro}
 |\Sha(E/\rat)| = O( \F_{E/\rat}^{7/4 + \epsilon(\F_{E/\rat})}),
 \end{equation} 
{\it where $\epsilon(\F)$ tends to $0$ when $\F$ tends to infinity.}

We give here bounds in these three directions, that is, for the product $|\Sha(A/K)| \Reg (A/K)$, the generators $\{P_i\}$ of the Mordell-Weil group, and the group $|\Sha(A/K)|$. Ours bounds are deduced from the following assumptions. 

\begin{hyp} \label{hyp-general} Let $A$ be an abelian variety of dimension $g$ defined over a number field $K$.   Suppose that $A$ carries a principal polarisation $\pp$. 
Suppose that the Tate-Shafarevic group $\Sha (A/K)$ is finite, the $L$-series of $A/K$ satisfies a condition on the order of growth (Hypothesis \ref{hyp-phragmen-lindeloef}), a functional equation (Conjecture \ref{funct-eq}) and the BSD-conjecture (Conjecture \ref{bsd}).
\end{hyp}

 \begin{rque}{\bf on the hypothesis.} The assumption that $A$ carries a principal  polarisation appears in Section \ref{section-local-periods} (in Lemma \ref{local-falt} where it is purely technical and unnecessary, and further in Lemma \ref{matrix-lemma}, to make evident the $\Im \tau_v$ and link the local periods appearing in the BSD-formula with the Faltings' height). Using Zarhin's trick, this hypothesis could be removed, up to modifying the constants in the last statements.  For this, one considers the variety $A^4 \times \Ad^4$, which carries a principal polarisation. See, e.g. \cite[Chap. XI (11.29)]{Edixhoven-Moonen-Van-der-Geer} for an explicit construction. 

Hypothesis \ref{hyp-phragmen-lindeloef} is less stringent than being of finite order, which  is satisfied in the modular case. In fact, in each case that Conjectures \ref{funct-eq} and \ref{bsd} are proven,  Hypothesis \ref{hyp-phragmen-lindeloef} is also proven (even if it doesn't formally  follow from Conjectures \ref{funct-eq} and \ref{bsd}). See Remark \ref{remark-hyp-phragmen-lindeloef}.
which however is satisfied in the modular case. (See Remark \ref{remark-hyp-phragmen-lindeloef}.)

\end{rque}
\smallskip

The core of our work is the next result, which goes in direction of Lang's conjecture (\ref{lang'sconj-Sha.Reg}). The bound of our proposition refines the conjectural bound (\ref{hindry'sconj}) of Hindry and extends Rémond's result \cite[Proposition A.2.3, Annexe A]{these.remond}, valid for $g=1$ and $K=\rat$. In what follows, we will state our results in a simplified form, which holds for $\F \ne 1$. See Section \ref{section-bounds-product-Sha.Reg} for more detailed bounds.

\begin{prop}\label{prop-sha.reg}
Suppose that $A/K$ satisfies Hypothesis \ref{hyp-general}. Then, with the above notations, when $\F \ne 1$,
\begin{equation} \label{borne-sha.reg}
|\Sha(A/K)| \times  \Reg(A(K) \leq (2^{16} g^2 d)^{\frac{gd}{2}} \cdot 2^{r} \cdot D_{K}^{g} \cdot \F^{\frac{1}{4} + \epsilon(\F)}  \cdot  e^{dh} \cdot \max\{1,h\}^{\frac{dg}{2}},
\end{equation}
where $\epsilon(\F) = 4gd \frac{\log \log \F}{\log \F} + 2gd \frac{\log \log (\F \cdot D_K^{2g})}{\log \F}$.  
\end{prop}

Using classical results on geometry of numbers on the euclidean structure provided by the  Néron-Tate height, and lower bounds for non-torsion points, we deduce our other results. 

 On one hand, we deduce a conditional upper bound for the N\'eron-Tate height of  the elements of a suitable basis of the Mordell-Weil group $A(K)$ modulo torsion.

\begin{thm} \label{thm-generateurs}
Suppose that $A/K$  satisfies Hypothesis \ref{hyp-general}. 
Then we can choose a system $\{P_{1}, \ldots, P_{r}\}$ of generators for the torsion-free part of the Mordell-Weil group $A(K)$ such that $\hat{h}_{\pp}(P_{1}) \leq \ldots \leq \hat{h}_{\pp}(P_{r}) $ and,  when $\F \ne 1$, 
\begin{equation}\label{borne-generateurs}
\hat{h}_{\pp}(P_r) \leq 
(2^{16} g^2 d)^{\frac{dg}{2} }\cdot (r!)^4  \cdot    D_{K}^{g}  \cdot  \F^{\frac{1}{4} + \epsilon(\F)}  \cdot 
e^{d h}\cdot \max \{d + g^g, h\}^{6074g(r-1) + \frac{dg}{2}},
\end{equation}
where $\epsilon(\F) = 4gd \frac{\log \log \F}{\log \F} + 2gd \frac{\log \log (\F \cdot D_K^{2g})}{\log \F}$. 
\end{thm}

See Section \ref{section-bounds} for more detailed bounds, and Remark \ref{comparaison-lang} for a comparison of our bound with Lang's conjecture (\ref{lang'sconj}) when $g=1$ and $K=\rat$.

On the other hand, we extend  Theorem 1 of \cite{goldfeld-szpiro}, the formula is (\ref{thm1-Goldfeld-Szpiro}),   to 
  principally polarised abelian varieties of arbitrary dimension, defined over an arbitrary number field. 
  When the dimension equals 1 and the number field is $\rat$, our bound improves (\ref{thm1-Goldfeld-Szpiro}).  
  See Remark \ref{comparaison-goldfeld-szpiro}, and, for detailed bounds, see Proposition \ref{BSD+Sz-Sha-detailed}.

\begin{thm}\label{BSD+Sz-Sha} Suppose that $A/K$ satisfies Hypothesis \ref{hyp-general}. 
Furthermore, suppose that $A/K$ satisfies  Szpiro's Conjecture (Conjecture \ref{gral.szpiro}). Then, for every $\epsilon >0$,
$$
|\Sha(A/K)| = O((\F_{A/K})^{\frac{1}{4} + \frac{1}{2} dg + \epsilon +  \delta(\F_{A/K})}),
$$
where $\epsilon >0$ and the implicit constant in the $O$ depends on $\epsilon$, $g$, $K$ and $r$, and $\delta(\F)$ tends to 0 when $\F$ tends to $\infty$.

\end{thm}

Ours results are stated with the Faltings' height of $A/K$, instead of the stable one, $h_{Stab}(A)$. In fact, when estimating the local periods, this height appears naturally (see Section \ref{section-local-periods}).  However, if the number field $K$ is large enough to contain the points of order 12, that is $A[12] \subset A(K)$, then we have equality, $h_{Stab}(A) = h_{Falt}(A/K)$ (see \cite[Exposé n.1, Corollaire 5.18]{sem-szpiro-pinceaux-ast86}).

In \cite{manin}, \cite{lang.conj.dio}, and \cite{goldfeld-szpiro} the argument is developed when the dimension is 1 and the number field is $\rat$, while  in \cite{hindry.mordell-weil} the dependence of the bounds on the number field is not the main interest and is not always made explicit. As pointed out in our joint work with Bosser \cite{bosser-surroca-brazil}, this dependence could play an important role.  For example, the discriminant of the number field appears in the rank of the variety. In fact, following the proof of the weak Mordell-Weil theorem, the latter can be bounded in terms of the logarithm of the discriminant of $K$ (\cite{ooe-top}).
Therefore, we consider here abelian varieties of arbitrary dimension over  arbitrary number fields, and make this dependence explicit. 
Furthermore, contrary to \cite{lang.conj.dio}  and \cite{hindry.mordell-weil} the bounds given here are not conjectured, but implied, by strong but, henceforth, classical conjectures.

As said before, the method is an extension of the one proposed by Manin, based on the BSD-conjecture. The BSD-conjecture  predicts the behavior of the $L$-series of the abelian variety $A$ at the center of symmetry, that is 1. In fact, it states that the order of vanishing of $L(A/K, s)$ at $s=1$ equals the Mordell-Weil rank of $A/K$. Furthermore, it gives a formula, which relates the value of the leading coefficient of the Taylor expansion of $L(A/K, s)$ at $s=1$ to the product of the Tate-Shafarevich group, the canonical regulator  and some other arithmetic invariants of the variety.  
The notations and the data concerning the abelian variety can be found in the next section. 
The core of our results are in Section \ref{section-lemmes}  where we bound the product of the Tate-Shafarevich group and the canonical regulator (Proposition \ref{prop-sha.reg}). 
In order to achieve this bound,  we estimate the remaining terms of the BSD-formula.
 To  deal with the leading coefficient of the Taylor expansion of the $L$-series  we use the functional equation (Lemma \ref{coeff-dominant}). We then relate  the local periods to the Faltings' height of $A/K$ (Lemma \ref{lemme-local-falt-ineg}).
We also give a bound for the torsion part of the Mordell-Weil group (Lemma \ref{lemma-torsion}). In Section \ref{section-lower-non-torsion} we recall some classical results on the geometry of numbers, and  comment on lower bounds for the N\'eron-Tate height of non-torsion points.
 In Section \ref{section-bounds} we deduce from the BSD-conjecture the bounds for the highest N\'eron-Tate height of a set of generators for $A(K)/A(K)_{\tors}$ (Theorem \ref{thm-generateurs}) and an upper bound for the order of $\Sha(A(K))$ (Theorem \ref{BSD+Sz-Sha}). We also give particular bounds in the one-dimensional case.

In fact, we apply these results in \cite{bosser-surroca-brazil} to an elliptic curve to show that, using the elliptic analogue of Baker's method in linear forms in logarithms, the BSD-conjecture for any single elliptic curve implies an inequality in the direction of the $abc$-conjecture over number fields.

\bigskip

Let's give a remark on the history of this paper. To deduce Theorem \ref{thm-generateurs} and Theorem \ref {BSD+Sz-Sha}  from Proposition \ref{prop-sha.reg}, we use a lower bound for the non-torsion points of the variety (see Section \ref{section-lower-non-torsion}).
A first version of this work, appearing in Axiv in 2008, uses a Masser's lower bound for the non-torsion points on a family of abelian varieties. Since the use of this bound was not appropriate in this context, we wrote a second version, using David's lower bound, together with the isogenies theorems of Masser-Whüstholz, for a simple principally polarised abelian variety. We then replaced it with a third version, using new isogenies theorems by Gaudron-Rémond. Finally, the lower bound used in this fourth version is the one by Bosser-Gaudron, which is completely explicit, and avoids the use of the isogenies theorems, simplifying the exposition. Some other sections where substantially modified, mostly  Sections \ref{section-bound-leading} and \ref{section-local-periods}.


\section{Notations}\label{notations}

Throughout  the text, we will consider an abelian variety $A$ of dimension $g$ defined over a number field $K$. We denote $d = [K:\rat]$ the degree and $D_K$ the absolute value of the discriminant of the field $K$. To $A$ one can  associate  different objects: the canonical regulator, the Tate-Shafarevich group, the conductor, the $L$-function, the Faltings' height.  
\medskip

\noindent {\bf 2.1.  The Néron-Tate height and the regulators.}

Let's denote $\Ad$ the {\it dual abelian variety} of $A$, that is, the connected component of the Picard group of $A$, denoted by $\Pic^{0}(A)$, which is also defined over $K$, and isogenous to $A$ (see [Mumford1970].) 
In particular 
\begin{equation}\label{duale}
 \dim (A) = \dim (\Ad), \hspace{0,15cm}
 \rk(A(K)) = \rk(\Ad (K)), \hspace{0,15cm}
 \F_{A/K} = \F_{\Ad/K}, \hspace{0,15cm}
 h_{\mathrm{Falt}}(A/K) = h_{\mathrm{Falt}}(\Ad/K).
\end{equation}
For the two last equalities see, respectively, \cite[Corollary 2]{serre-tate} and, \cite[Corollaire 2.1.3]{raynaud-asterisque127-1985}. 

\medskip

Let $\cal L$ be an invertible sheaf on $A$ and let $\tau_x : A \to A$ be the translation by $x \in A$. We define a morphism $\phi_{\cal L} : A \to \Ad$ by $\phi_{\cal L} (x) = \tau_x^* {\cal L} \otimes {\cal{L}}^{-1}$. The morphism $\phi_{\mathcal{L}}$ is an isogeny  (i.e., is surjective and has finite kernel) if and only if $\cal{L}$ is ample. In this case, we call it a {\it polarisation} and  $\deg(\phi_{\cal L}) := |\ker \phi_{\mathcal{L}}| = h^0(A, {\cal L})^{2}$. The polarisation is called {\it principal} if its degree is 1. 

\medskip

 Let us fix such an ample line bundle $\pp$ on $A$, which is identified with an ample symmetric divisor.  
 The associated Néron-Tate height on $A(\Kbarre)$ (see, e.g. \cite[Section B.5.]{hindry-silverman}), 
 $$\hat{h}_{\pp}  := \hat{h}_{A,\pp} : A(\Kbarre) \to \R,$$ 
is a positive quadratic form. Since $\hat{h}_{\pp}(P) = 0$ if and only if $P$ is a torsion point, $\hat{h}_{\pp}$ is a positive definite quadratic form on $A(K)/A(K)_{\tors}$, and one could prove that $\hat{h}_{\pp}$ 
extends to a positive definite quadratic form on  $A(K) \otimes_{\Z} \R$. The associated bilinear pairing 
$$<P,P'>_{\pp} := \frac{1}{2} (\hat{h}_{\pp}(P+P') - \hat{h}_{\pp}(P) - \hat{h}_{\pp}(P')),$$ 
which satisfies $<P,P>_{\pp} = \hat{h}_{\pp} (P)$,  endows $A(K) \otimes_{\Z} \R \simeq \R^{r}$   with a structure of a euclidean space, and we can view $A(K)/A(K)_{\tors}$ as a lattice sitting inside this space. The regulator 
 $$\Reg_{\pp}(A/K) :=\det (<P_i,P_j>_{\pp} )_{1\leq i,j \leq r} \geq 0,$$
where $\{P_1, \ldots, P_r\}$ is a basis  for $A(K)/A(K)_{\tors}$,  is the square of the volume of the fundamental domain for the lattice. 

\medskip

 We will define now a {\it canonical} regulator, independent of the choice of the line bundle $\pp$. 
 
 Let $\cal P$ be the Poincar\'e line bundle on $A \times \Ad$ and $\hat h_{\cal P} : = \hat{h}_{A\times \Ad,\mathcal{P}}$ the canonical height on $A \times \Ad$ with respect to $\cal P$. Let's define a pairing, for $P \in A(K)$ and $Q\in \Ad(K)$,
 $$<P,Q>_{\cal P} := \hat h_{\cal P}(P,Q).$$
  Choose a $\enteros$-basis $\{P_1, \ldots, P_r\}$, resp. $\{Q_1, \ldots, Q_r\}$,   for  $A(K)/A(K)_{\tors}$, resp. of $\Ad(K)/\Ad(K)_{\tors}$. The {\it canonical regulator} of $A$  (also called the {\it discriminant of the height pairing}) is defined by
 $$
\Reg(A) : = \left |  \det((<P_{i}, Q_{j}>_{\cal P})_{1\leq i, j\leq r})\right |.
 $$
It is a non-zero real number and does not depend on the choice of the basis (see \cite{tate.bourbaki.bsd}). However, it could be related to $\Reg_{\pp}(A/K)$. In fact, we could recover all heights on $A$ from the canonical pairing  (see, e.g., \cite[Remark p.39]{serre} or \cite{hindry.mordell-weil}): 
\begin{equation}\label{all-heights}
 \hat h_{\pp} (P) = -\frac{1}{2} <P,\phi_{\pp}(P)>_{\cal P}. 
  \end{equation} 
Denotes   $$u = [\Ad(K)/\Ad(K)_{\tors}:\phi_{\pp}(A(K)/A(K)_{\tors})] $$ 
the index of the subgroup $\phi_{\pp}(A(K)/A(K)_{\tors})$ in $\Ad(K)/\Ad(K)_{\tors}$. The index $u$ equals 1 if the polarisation $\pp$ is principal, and we can prove that, in general, $u$ divides $\deg(\phi_{\pp} )^r$. Using $(\ref{all-heights})$ we can prove that
 \begin{equation}\label{reg}
 \Reg_{\pp} (A/K) = u 2^{-r} \Reg (A/K)  \leq 2^{-r}   \deg(\phi_{\pp})^r \Reg(A/K).
 \end{equation}

  \medskip

\noindent {\bf 2.2.  The Tate-Shafarevich group.} The {\it Tate-Shafarevich group} of $A/K$ is defined by
$$\Sha (A/K) : = \ker (H^1(\Gal(\Kbarre/K), A_K) \rightarrow \prod_v H^1(\Gal(\Kbarre_v/K_v), A_{K_v})).$$ 
Rubin  \cite{rubin1987} gave the first examples of elliptic curves for which it can be proved that the Tate-Shafarevich group is finite  for elliptic curves defined over $\rat$. See also the results of Kolyvagin  \cite{kolyvagin1989}.  

We  suppose throughout that $\Sha (A/K)$ is finite.

\medskip 

\noindent {\bf 2.3. The conductor.} 
Let $v$ be a finite place of $K$, which corresponds  to a prime ideal $\vp$. We will denote $K_v$ or $K_{\vp}$ the completion of $K$ at $v$. 
 For any prime ideal $\vp$ of $K$, fixing a prime in $K_{\vp}$ above $\vp$, gives us a decomposition group $G_{\vp} = \Gal (\overline{K_{\vp}}/K_{\vp})$ for $\vp$ in $\Gal(\Kbarre/K)$. Let $I_{\vp}$ be the inertia subgroup of $G_{\vp}$, inducing the identity on the residue field $k(\vp)$. Let $\pi_{\vp}$ denote the Frobenius, which generates the quotient $G_{\vp}/I_{\vp}$. Up to conjugation, $G_{\vp}, I_{\vp}$ and $\pi_{\vp}$ depend only on $\vp$. 
 Let $\ell$ be any prime,  $\ell \ne char ( k(\vp) )$. Denote $A[N]$ the $N$-torsion of $A$, for an integer $N$, $T_{\ell}(A)= \varprojlim A[\ell^n]$ the $\ell$-adic Tate module, and $V_{\ell}(A/K)= T_{\ell}(A)\otimes_{\Z_{\ell}}\rat_{\ell}$ the $\rat_{\ell}$-vector space associated.  Since $\Gal(\Kbarre/K)$ acts on  $V_{\ell}(A/K)$,  we have a $\ell$-adic representation $\rho : \Gal (\Kbarre/K) \rightarrow GL(V_{\ell}(A/K))$.

The {\it conductor} of the abelian variety $A/K$ is the integral ideal of $K$ defined by
 $$\F_{A/K} = \prod \vp^{f_{\vp}},$$
 where the product runs over the prime ideals $\vp$ of $K$, and $f_{\vp}$ is a positive  integer, called the {\it exponent of the conductor},  which we will define below. The exponent $f_{\vp}$ is zero if and only if $A$ has good reduction at $\vp$.   As in  \cite{serre.zeta}, we will attach to the representation $\rho$ two positive integers $\varepsilon_{\vp}(\ell)$ and $\delta_{\vp}(\ell)$, which measure the ramification of $\rho$. We follow the notations of \cite{lockhart-rosen-silverman}. Denote  $V_{\ell}(A/K)^{I_{\vp}}$the submodule of elements fixed by $I_{\vp}$.  Define
 $$\varepsilon_{\vp}(\ell)= \codim_{\rat_{\ell}}V_{\ell}(A/K)^{I_{\vp}}.$$
 Let $L_{\vp}= K_{\vp}(A[l])$ be the field generated over  $K_{\vp}$  by the $\ell$-torsion points of $A$. Denote $v_{L_{\vp}}$ the normalised valuation on $L_{\vp}$. Let $\pi_{L_{\vp}}$ be a uniformiser for $L_{\vp}$. Denote $G_i= \{\sigma \in \Gal(L_{\vp}/K_{\vp}); v_{L_{\vp}}(\sigma \pi_{L_{\vp}} - \pi_{L_{\vp}}) \geq i+1 \}$ the $i$-th inertia group associated to $L_{\vp}/K_{\vp}$  and $g_i = |G_i|$ its order. Write $g_0 = |\Gal(L_{\vp}/K_{\vp})|$. We define
 $$\delta_{\vp}(\ell)= \sum_{i\geq 1} \frac{g_i}{g_0} \dim_{\mathbf{F}_l}\left(\frac{A[l]}{A[l]^{G_i}}\right).$$
 It has been proven (see the references in \cite{lockhart-rosen-silverman}) that $\varepsilon_{\vp}(\ell)$ and $\delta_{\vp}(\ell)$ are independent of $\ell$  so we will denote them by $\varepsilon_{\vp} $ and $\delta_{\vp}$. They are called the {\it tame} part and the {\it wild} part of the conductor, respectively. The exponent of the conductor is given by
 $$f_{\vp}= \varepsilon_{\vp} + \delta_{\vp}.$$

   It is known that if $p> 2g+1$, where $p$ is the prime number lying below $\vp$, then $f_{\vp} \leq 2g$. Furthermore, it is proven in \cite{lockhart-rosen-silverman} that $f_{\vp} \leq 12 g^2 v_{K_{\vp}} (p)$ unconditionally (see \cite{brumer-kramer1994} for the best possible upper bounds in all cases).
   
 It is also known that for elliptic curves defined over $\rat$, the conducteur satisfies $\F_{E/\rat} \geq 11$. 
In higher dimension, we still have the following  lower bound over $\rat$ : $\F_{A/\rat} > 3$ (see {\it loc. cit.}).  In revanche, there are number fields where the abelian variety have good reduction everywhere (see \cite{schoof-2003}) and then, $\F_{A/K} = 1$.  

  \medskip
  
 \noindent {\bf 2.4. $L$-series.}
 We now define  the $L${\it -series}, also called the $\zeta${\it -function}, of the variety $A$ (see \cite[Section 4]{serre.zeta}). 
 Since the Frobenius is defined up to $I_{\vp}$, it makes sense to define a polynomial $P_{A, \vp} (T) = \det (1 - (\rho (\pi_{\vp})  | V_{\ell}(A/K)^{I_{\vp}})T)$, where $\pi_{\vp}$ is regarded as acting on the submodule $V_{\ell}(A/K)^{I_{\vp}}$ of elements fixed by $I_{\vp}$. The polynomial $P_{A, \vp} (T)$ has integral coefficients which are independent of $\ell$ (\cite[Theorem 3]{serre-tate}). 
Define
$$L(A/K,s) = \prod _{v_{\vp}} P_{A,\vp}(N(v_{\vp})^{-s})^{-1}$$ 
where the product is taken over all non-archimedean places $v_{\vp}$ of $K$ and $N(v_{\vp})$ is the norm of the prime ideal $\vp$ associated to $v_{\vp}$. Define the {\it normalised} $\ell${\it -function} by %
$$\Lambda(A/K,s) = (N_{K/\rat}( \mathcal{F}_{A/K}) \cdot D_{K}^{2g})^{s/2}  \cdot ((2\pi)^{-s} \cdot \Gamma (s))^{g[K:\rat]} \cdot L(A/K, s),$$
where $\Gamma(z) = \int_0^{+ \infty} t^{z-1} e^{-t} dt$ is the classical $\Gamma$-function. For the $\Gamma$-factors see \cite[section 3]{serre.zeta}.
  The Euler product converges and gives an analytic function for all $s$ satisfying $\Re(s) > \frac{3}{2}$. 
  
  \medskip

\noindent{\bf 2.5. The local factor.}
In order to state the Birch and Swinnerton-Dyer conjecture, we will introduce a {\it local factor} for $A$. 
We follow here the approach of Gross \cite{gross-BSD82}, using the theory of N\'eron models and Tamawaga numbers.\footnote{For a general formulation on BSD, we could see the original formulation of Tate on abelian varieties \cite{tate.bourbaki.bsd}, and also \cite{milne}, as well as \cite{bloch1980} for a volume theoretic formulation.}

Let $\A$ denote the N\'eron model of $A$ over the ring of integers $\ok$ of $K$ and $\A^0$ the largest open subgroup of $\A$ in which all fibers are connected. 
 
 As a first step, to every place $v$ of $K$, we will associate a local number  $c_v$. 
 
 For a {\it finite} place $v$,  $\A(K_v)$ is a commutative group, and  $\A^0(K_v)$ is the subgroup of the $K_v$-rational points which reduces to the identity component of the N\'eron model $\A$. %
Denote
\begin{equation}\label{finite-local-factor}
c_v:=  \left| \A(K_v):\A^0(K_v) \right|
\end{equation}
the index. 
Since this integer is 1 for almost all $v$ (that is, for the places where $A$ has good reduction), we may define the product
\begin{equation}c_f(A/K) = \prod_{v \in M_K^0} c_v.
\end{equation}

Denote $\DiffAK$ the sheaf of differential 1-forms on $A/K$.   Let  $\{\omega_{1}, \ldots, \omega_{g}\}$ be a $K$-basis of $H^{0}(A, \DiffAK)$.  Then  $\eta= \omega_{1}\wedge \ldots \wedge \omega_{g}$ is a non-zero differential $g$-form on $A$. Let $\DiffA$ be the invertible sheaf of the differential $g$-forms on $\A$.

 The module $H^{0}(\A, \DiffA)$ of global invariant differentials on  $\A$ is  a projective $\ok$-module of rank 1  and can be written as
$$H^{0}(\A, \DiffA) = \eta \cdot \va,$$
where $\va = \va_{\eta}$ is a fractional ideal of $K$ (depending on $\eta$). We have 
$$\va_{\eta} = \va_{\alpha \eta} \cdot (\alpha)\esp \textrm{ for } \esp \alpha \in K^{\star}.$$

Let $v$ be an {\it archimedean} place of $K$ and let's denote $A_v$ the abelian variety obtained  from $A$ by the action by $v$.  The integral homology $H_1(A_v(\overline{K_{v}}), \Z)$ of $A_v$  is a free $\enteros$-module of rank $2g$. 
Let $(\gamma_{1,v}, \ldots, \gamma_{2g,v})$ be a basis of  $H_1(A_v(\overline{K_{v}}), \Z)$. 

Let $v$ be a {\it complex} place and define, $\omega_{g +j} = \overline{\omega_{j}}$ for $j\in \{1, \ldots, g\}$, and
$$c_v := c_v(A, \eta) = \left|\det \left(\int_{\gamma_{i,v}} \omega_{j} \right)_{1 \leq i, j \leq  2g}\right|.$$

Let $v$ be a {\it real} place. Choose the basis  so that $\gamma_{1,v}, \ldots, \gamma_{g,v}$ generates the part of $H_1(A_v(\overline{K_{v}}), \Z)$ fixed by complex conjugation. 
Denote  $|A_{v}(\R) : A_{v}(\R)^{0}|$ the number of real components of the variety $A_{v}$ and define$$c_v := c_v(A, \eta) = |A_{v}(\R) : A_{v}(\R)^{0}| \cdot \left|\det \left (\int_{\gamma_{i,v}} \omega_{j} \right)_{1 \leq i, j \leq g}\right|.$$

These two integrals are non-zero and depend only on $A$, the $g$-form $\eta$ and the place $v$. 

 We also define  the {\it archimedean local factor} as 
\begin{equation}\label{archimedean-local-factor}
c_{\infty}(A/K) = N_{K/\rat}(\va) \cdot \prod_{v \in M_K^{\infty}}c_v,  
\end{equation}
which is independent of the choice of the differential $\eta$ (because of the product formula). Finally, we define the {\it local factor} of $A$ as 
\begin{equation}
c(A/K) = c_{f}(A/K) c_{\infty}(A/K),
\end{equation}
which is a positive real number.

  Manin \cite{manin} gives another formulation for the local factors (for $g=1$), in terms of measures. 
We follow here \cite[Chapitre II]{these.remond} (see pages 15 and 18, where he denotes $\omega = \omega_1 \wedge \ldots \wedge \omega_g$ our $\eta$, and  $m_v(A,\omega)$ our $c_v = c_v (A, \eta))$. See also \cite[p.224]{gross-BSD82}) and \cite{tate.bourbaki.bsd}.  

Indeed, let $\mu_v$ be an additive Haar measure on $K_v$ such that $\mu_v(O_{K_v}) = 1$ if $v$ is finite, $\mu_v$ is the Lebesgue measure if $v$ is a real archimedean place (i.e. $K_{v} = \R$) and twice the Lebesgue measure if $v$ is complex (i.e. $K_{v} = \C$). 
Then, when the differential form $\eta \in H^0(A, \Omega^g_{A/K})$ is algebraic, it induces an analytic form to which one could associate the {\it module measure} $\mod(\eta)$, denoted by $|\eta| \mu_v^g$ by Gross, and called the associated Haar measure on $A(K_v)$. We then have the following formulas relating the local periods $c_v(A, \eta)$ with the module measure  (see \cite[p.224]{gross-BSD82} and \cite{these.remond} for further details). For an archimedean place, we have (\cite[p.15]{these.remond})

\begin{equation}\label{mod-measure-archimedean}
c_v(A, \eta) = \int_{A(K_v)} \mod(\eta),
\end{equation}
and also the formula (see \cite[p.18]{these.remond})
\begin{equation}\label{mod-measure-wedge-complex}
\int_{A(\C)} \mod(\eta) = \int_{A(\C)} |\eta \wedge \overline{\eta}|.
\end{equation}
For a finite place, we have the formula
$$ \int_{A(K_v)} \mod(\eta) = N(v)^{\ord_v(\va_{\eta} )-g} \cdot |\A(K_v)|.$$
Since $P_v(N_{K/\rat}(\vp_v)^{-s})^{-1} = N(v)^g \cdot |\A^0(K_v)|^{-1}$, we then obtain, for $v$ finite,

$$c_v (A,\eta) =P_v(N_{K/\rat}(\vp_v)^{-s})^{-1} \cdot N(v)^{-\ord_v(\va_{\eta} )}   \int_{A(K_{v})}  \mod(\eta).$$

\medskip

\noindent {\bf 2.6. The Faltings' height.} The part of the BSD-formula concerning the local periods $c_v$ can be bounded in terms of the {\it Faltings' height}. We define it in the following way, using Faltings' normalisation, $\frac{1}{2^g}$, as in  \cite[Chapter II]{cornell-silverman} (where we corrected the typo).

We endowed the line bundle $\DiffA$ with an Hermitian metric by defining, for a section $s$ and for every archimedean place $v$,
\begin{equation}\label{norme2}
|s|_{v} := \left( \frac{1}{2^g} \int _{A(\overline{K_{v}})} \left | s \wedge \overline{s}   \right | \right)^{\frac{1}{2}}.
\end{equation}
We also define 
$$|| s ||_{v} = |s|_{v}^{n_{v}},$$ 
where $n_{v} = [K_v:\rat_v]$ equals 1 if $v$ is real and   equals 2 if $v$ is complex.
We remark that Rémond \cite[p.17]{these.remond} have different notations and normalisation than we use here. In fact, he denotes $||s||_v$ the integral $\int_{A(\C)} |s \wedge \overline{s} |$ which corresponds to our $2^g \cdot |s|_v^2$.
However, this norm extends the norm on $K_{v}$  (i.e. $\forall k \in K_{v}, \esp \forall s \in \DiffA \otimes_{\ok} K_{v}, \esp ||ks ||_{v} = ||k||_{v} \cdot ||s||_{v}$).

Taking the pull-back of $\DiffA$ and metrics via the neutral section $e : \Sok \rightarrow \A$, we obtain a metrised line bundle on $\Sok$ (i.e. a projective $\ok$-module of rank 1): 
$$\eDiffA := e^{*}\DiffA.$$

The line bundle $\eDiffA $ can be identified with $H^{0}(\DiffA) = \eta \cdot \va_{\eta}$. In fact, $\eDiffA = e^{*}\DiffA = \pi_{*} \DiffA$, where $\pi : \A \rightarrow \Sok$ is the structural morphism, and since the line bundle is affine, it can be identified with  the module of its global sections $H^{0}(\DiffA)$.  %

The {\it Faltings' height} of $A$ is the {\it Arakelov degree} of $\omega_{\A/\ok}$ :
$$h_{\mathrm{Falt}}(A/K) = \frac{1}{[K:\rat]} \deg_{\mathrm{Ar}}(\eDiffA, ||.||) = - \frac{1}{[K:\rat]} \log \prod_{v \in M_{K}} ||s||_{v},$$
for any section $s$. 
The Faltings' height denoted by Rémond by $h(A)$ corresponds, with our notations, to $h_{Falt}(A/K)
 - \log(2^{g/2})$. 
 
 It is well known that %
$$\deg_{\mathrm{Ar}}(\eDiffA, ||.||) = \log \card (\eDiffA/s\ok) - \sum_{v | \infty} \log ||s||_{v}.$$
%

The height defined in the same way but, over a number field extension where $A$ has semi-stable reduction, is called {\it stable} Faltings' height. We denote it by $h_{stab}(A)$. It doesn't depends on the ground field and it satisfies
$$h_{stab}(A) \leq h_{Falt}(A/K),$$
with equality if, and only if, $A/K$ is semi-stable.

\section{On the Birch and Swinnerton-Dyer conjecture}\label{section-lemmes}

We can now give a classical generalisation of a conjecture of Hasse-Weil and state the celebrated conjecture of Birch and Swinnerton-Dyer  (see, e.g.,  \cite{birch-swinnerton-dyer} for the case of elliptic curves and \cite{gross-BSD82} for a general formulation for an abelian variety defined over an arbitrary number field).


\begin{conj} [Hasse-Weil] \label{funct-eq} 
Let $A/K$ be an abelian variety defined over a number field. The $L$-series and the $\Lambda$-series of $A/K$ have an analytic continuation  to the entire complex plane and the $\Lambda$-series satisfies the functional equation 
$$\Lambda(A/K, 2-s) = \varepsilon \Lambda (A/K, s), \esp \textrm{for some }\esp \varepsilon = \pm 1.$$
\end{conj}

This conjecture is true for abelian varieties with complex multiplication (\cite{shimura-taniyama}), in some special cases, this conjecture is also true for modular abelian varieties (\cite{shimura.automorphic.book}) and it is true for elliptic curves over $\rat$ (\cite{wiles1995} and \cite{breuil-conrad-diamond-taylor}). See also \cite{patrikis-taylor2015}  and  \cite{BCGP2018} and, for elliptic curves defined over a real quadratic  field, \cite{freitas-bao-siksek2015}.

\begin{conj}[Birch and Swinnerton-Dyer]\label{bsd} Let $A$ be an abelian variety defined over a number field $K$.
\begin{enumerate}
\item The $L$-series $L(A/K,s)$ has an analytic continuation to the entire complex plane.
\item $\ord_{s=1}L(A/K,s) = \rk (A(K))$.
\item The leading coefficient $L^{\star}(A/K, 1) = \lim_{s\to1}\frac{L(A/K,s)}{(s-1)^{\rk(A(K))}}$ in the Taylor expansion of $L(A/K,s)$ at $s=1$ satisfies 
\begin{equation}\label{formula-bsd}
 L^{\star}(A/K, 1) =   | \Sha (A/K)| \cdot \Reg(A(K)) \cdot |A(K)_{\tors}|^{-1} \cdot |\Ad(K)_{\tors}|^{-1} \cdot c(A/K) \cdot D_K^{-g/2}.
\end{equation}
\end{enumerate}
\end{conj}

 In the 70's and 80's, striking progress was achieved on the BSD-conjecture, providing evidence for its truth (\cite{coates-wiles1977}, \cite{gross-zagier1986}, \cite{rubin1987}, \cite{kolyvagin1989}). In particular, for an elliptic curve defined over $\rat$ satisfying $\ord_{s=1}L(E/\rat,s) = 0$, the conditions 1. and 2. are proved and also a relation between the value of $L(E/\rat,1)$ and the order of $\Sha(E/\rat)$ similar to condition 3., up to some factor term. More recently,  more evidence arises from \cite{Bhargava-Shankar2015}. See also the references in \cite{kato-trihan2003}, in particular for the function field case, where much more is known. 

 In this section we bound the product $|\Sha(A/K)| \cdot \Reg(A/K)$ from above. In order to do it, the formula (\ref{formula-bsd}) of the BSD-conjecture suggests to bound each one of the remaining terms. This is done in the following  lemmas.

\subsection{Bound for the leading coefficient $L^{*}(A/K, 1)$}\label{section-bound-leading}

For his algorithm,  Manin  deals with the case when $A=E$ is an elliptic curve and $K=\rat$ (\cite[Theorem 11.1]{manin}). Following his notations, let's consider the Dirichlet expansion of the $L$-series  $L(E/\rat, s) = \sum_{n=1}^{\infty}a_n n^{-s}$ and set $F(z) =  \sum_{n=1}^{\infty} a_n e^{2\pi i n z}$. Then $F$ is holomorphic in the upper half plane. He then uses (several times) the Hecke functional equation $F(z) = \varepsilon N^{-1} z^{-2} F(-\frac{1}{Nz})$, where $N >0$ is the conductor of $E$, and the fact that the sequence $(a_n)$ does not grow faster than $O(n^c)$ for some $c> 0$. (In fact, in this case, $|a_n| \leq n^{1/2} \tau(n)$, where $\tau(n)$ is the number of divisors of $n$.) See also Theorem 9.3 a), Remarks 9.7 and Section 11.3 of \textit{loc. cit.}.

Classically, when looking forward to bound the leading coefficient of such $L$-series, one proceeds in three steps. First, a bound could be easily given for $| \Lambda (s)|$ in the half-plane defined by $\Re(s) > 3/2$, using the Hasse-Weil bound. Second, using the functional equation, one proves that the bound is still valid in the other half-plane defined by $\Re(s) < 1/2$. What is missing, is to bound it in the vertical strip defined by $1/2 < \Re(s) < 3/2$. However, in each case when Conjectures \ref{funct-eq} and \ref{bsd} are proven,  it is also proven that $|\Lambda(s)|$ is bounded in this vertical strip (even if it doesn't follow formally from Conjectures \ref{funct-eq} and \ref{bsd}, but, e.g. from modularity).

Here we use the Hasse-Weil bound for bounding $|\Lambda(s)|$ in the right-half-plane $\Re(s) > 3/2$, the functional equation for the $\Lambda$-series (Conjecture \ref{funct-eq}) to bound it in the left-half-plane $-1/2 < \Re(s)$, and  the classical convexity argument as in the Phragm\'en-Lindelöf principle (see \cite{phragmen-lindeloef1908} or \cite[section 5.61 page 177]{titchmarsh-book-th-funct}) to bound it in the vertical strip $-1/2 \leq \Re(s) \leq 3/2$.
Nevertheless, in order to apply  the Phragm\'en-Lindel\"of principle, we will use a  hypothesis on the order of growth of the $\Lambda$-series, Hypothesis \ref{hyp-phragmen-lindeloef}, as in (\ref{condition-phragmen-lindeloef-strip}) below. In some sense, this condition replaces the condition on the growth of the sequence $(a_n)$ of the coefficients of the $L$-series used by Manin. 
We conclude by applying the Cauchy inequality in two different ways so as to obtain different kinds of bounds for the leading coefficient $L^*(A/K,1)$.

\smallskip

\begin{lemma}[Phram\'en-Lindel\"of]\label{phragmen-lindeloef-lemma}

Let $f(z)$ be an analytic function of $z = r e^{i \theta}$, regular in the region $D$ between two straight lines making an angle $\pi/\alpha$ at the origin, and on the lines themselves. Suppose that 
$$f(z) \leq M$$
on the lines, and that, 
\begin{equation}\label{condition-phragmen-lindeloef}
f(z) = O (e^{r^{\beta}}), \textrm{as} \esp  r \to \infty, \textrm{where} \esp \beta < \alpha,
\end{equation}
uniformly in the angle. Then actually, the inequality $f(z) \leq M$ holds throughout the region\nolinebreak[4] $D$.
\end{lemma}

We will use the above result when the angle is transformed into a strip.

\begin{lemma}\label{phragmen-lindeloef-lemma-strip}

Let $\epsilon > 0$ be any positive number. Set $\alpha(\epsilon) = \frac{\pi}{2 \epsilon +1}$. Let $\phi(s)$ be an analytic function, regular in the strip $S$ between the two parallel lines $\sigma = 3/2 + \epsilon$ and $\sigma = 1/2 - \epsilon$, and on the lines themselves. Suppose that 
$$\phi(s) \leq M$$
on the lines, and that, 
\begin{equation}\label{condition-phragmen-lindeloef-strip}
\phi (s) = O(e^{e^{\rho \tau}}), \textrm{as} \esp  \tau = \Im(s) \to \infty, \textrm{where} \esp \rho < \alpha (\epsilon),
\end{equation}
uniformly in the angle. Then actually, the inequality $\phi(s) \leq M$ holds throughout the strip $S$.
\end{lemma}

\noindent {\it Proof of Lemma \ref{phragmen-lindeloef-lemma-strip}.} 
We use the notation $z=r e^{i \theta}$, $r \geq 0$, in the region $D$ of the complex plane defined by $|\theta| \leq \frac{\pi}{2\alpha(\epsilon)}$. We set $s= i \log z +1$, and $f(z) = \phi(s)$. In this way, the  two straight lines making an angle $\pi/2$ at the origin are transformed into two parallel lines $\sigma = 3/2 + \epsilon$ and $\sigma = 1/2 - \epsilon$. The origin is sent to the ``infinity" of the negative part of the imaginary axis. The imaginary part of our variable $\sigma$ is $\tau = \log r$. 
 Since $\phi(s)$ satisfies condition (\ref{condition-phragmen-lindeloef-strip}), the condition (\ref{condition-phragmen-lindeloef}) on the growth of the function $f(z)$ is satisfied. We conclude applying Lemma \ref{phragmen-lindeloef-lemma}. 
\hfill $\Box$

\begin{hyp}\label{hyp-phragmen-lindeloef}
Let $A/K$ be an abelian variety defined over a number field. The $\Lambda$-series of $A/K$ satisfies 
$$|\Lambda(A/K,s)| = O(e^{e^{\rho \tau}}), \textrm{as }  \tau = \Im(s) \to \infty,  \textrm{where } \rho \leq \frac{\pi}{2}, uniformly.$$
\end{hyp}

\begin{remark} \label{remark-hyp-phragmen-lindeloef}
{\em{
Hypothesis \ref{hyp-phragmen-lindeloef} is less stringent than being  \textit{of finite order}\footnote{The order of the function $f$ is $\inf_m \{f(z) = O(e^{|z|^m}), \textrm{as } |z| \to \infty \} = \lim \sup_{R \to \infty} \frac{\log \log \max_{|z| \leq R}|f(z)|}{\log R}$.}. One expects that the $\Lambda$-series of an abelian variety is always of finite order.  
 For example, $\Lambda$ is of order 1 in the case of an elliptic curve over $\rat$. Indeed, as soon as the $L$-series  $L(s) = \sum_{n=1}^{\infty}a_n n^{-s}$ is the Mellin transform of a function $F(z) =  \sum_{n=1}^{\infty} a_n e^{2\pi i n z}$,  which is modular, and satisfies $a_n = O(n^c)$, one can apply Hecke's lemma (see \cite[Lemma 1]{weil-math-ann1967}, B $\implies$ A) to bound  $|\Lambda(s)|$ in a vertical strip, say $-3 \leq \Re(s) \leq 5 $. We then bound $|\Lambda(s)|$ for $\Re(s) \geq 5$. Using again the Hecke functional equation for $F$, the bound also holds in the half plane $\Re(s) \leq -3$. Then,  for $|s| = R \gg1$,  $|\Lambda(s)| \ll e^{R\log R}$ and we conclude that  the normalised $\Lambda$-series is of order $1$.

As it is well known, the further condition (\ref{condition-phragmen-lindeloef-strip}), and thus Hypothesis \ref{hyp-phragmen-lindeloef}, is required to apply the Phragm\'en-Lindel\"of theorem. Indeed, wild functions such as $e^{e^z}$, on the strip $|\Im(z)| < \frac{\pi}{2}$, are not bounded by their maximum on the boundary. This function doesn't satisfy (\ref{condition-phragmen-lindeloef-strip}). }}
\end{remark}

\begin{lemma}\label{coeff-dominant}

Let $A/K$ be an abelian variety of dimension $g$ satisfying Conjecture \ref{funct-eq} and Hypothesis \ref{hyp-phragmen-lindeloef}.  Let $r= \ord_{s=1}L(A/K,s)$ and  $\mathcal{F} = N_{K/\rat}(\mathcal{F}_{A/K})$. 
Then the leading coefficient of the $L$-series of $A/K$ at $s=1$ satisfies the following bounds:
\begin{equation} \label{bound-coeff-rank}|L^{\star}(A/K, 1)| \leq \left(9/2\pi\right)^{g[K:\rat]}\sqrt{\mathcal{F}} \cdot D_{K}^{g}
\end{equation}
\begin{equation} \label{bound-coeff-cond}
|L^{\star}(A/K, 1)| \leq  e \cdot 2^{r}  \cdot (6/5)^{g[K:\rat]}\cdot \F^{\frac{1}{4}}\cdot D_{K}^{\frac{g}{2}}\cdot (\log (\F \cdot D_{K}^{2g}))^{2g[K:\rat]}.
\end{equation}
\end{lemma}

Before proving the lemma, let us explain why we give two different bounds.

\begin{remark} \label{remark-rank}
{\em{The upper bounds (\ref{bound-coeff-rank}) and (\ref{bound-coeff-cond}) depends on $g, [K:\rat], \F$ and $D_K$. The bound (\ref{bound-coeff-cond}) also depends on the order of the $L$-series at 1, here denoted by $r$. Conjecturally,  $r$ equals the rank of $A(K)$. 
When $g=1$ and $K=\rat$, it is expected that $\rk(E(\rat))=0$ or $1$. In this case, the bound (\ref{bound-coeff-rank}) is sharper than (\ref{bound-coeff-cond}) as soon as $\mathcal{F}\geq 6$. Moreover, 
as pointed out in \cite{bosser-surroca-brazil}, the dependence of these bounds on the number field could play an important role in some applications. 
Concerning the rank, Ooe and Top \cite{ooe-top} proved the following bound:
\begin{equation}\label{bound.ooe-top}
\rk(A(K)) \leq \gamma_1 \log \F + \gamma_2 \log D_K + \gamma_3,
\end{equation}
where $\gamma_1, \gamma_2$ and $\gamma_3$ are positive real numbers depending only on $g$ and $[K:\rat]$ (see \cite[Proposition 5.1]{remond-london2010} for explicit computations of $\gamma_1, \gamma_2$ and $\gamma_3$).  Using (\ref{bound.ooe-top}), we deduce from (\ref{bound-coeff-cond})  a bound independent of the rank, which growth in $\F$ and $D_K$ is as: 
$$\F^{\frac{5}{4}}\cdot D_K^{\frac{g}{2} +1}\cdot (\log (\F\cdot D_K^{2g}))^{2g[K:\rat]}.$$
With respect to the conductor, the estimate (\ref{bound-coeff-rank}) is of better quality than this last one. 
As for the dependence on the discriminant $D_K$,  the estimate (\ref{bound-coeff-rank}) also has a better dependence than this last bound if, and only if,  the dimension $g$ is 1 or 2. 

Depending on our focus, we will use (\ref{bound-coeff-rank}) or (\ref{bound-coeff-cond}).  E.g. in \cite{bosser-surroca-brazil} we are concerned with elliptic curves and we are interested on the dependence on $D_K$: the bound (\ref{bound-coeff-rank}) is used therein. However, it is expected that $\rk(A(K)) \ll \frac{\log \F}{\log \log \F}$. This would give, using (\ref{bound-coeff-cond}), a bound for the leading coefficient of the order
$$\F^{\frac{1}{4} + \epsilon(\F)}\cdot D_K^{\frac{g}{2} + 1 + \epsilon'(D_K)},$$
where $\epsilon$ and $\epsilon'$ depend on $g$ and $[K:\rat]$ and tend to $0$ when $\F $ tends to infinity and, respectively, when $D_K$ tends to infinity, which  is sharper than (\ref{bound-coeff-rank}).
}}
\end{remark}

\medskip

\noindent {\it Proof of Lemma \ref{coeff-dominant}.} 
Let us consider  the abelian variety $A' = \Res^{K}_{\rat} A$ over $\rat$, which is obtained from $A$ by restriction of scalars (see \cite{milne}).  Over $\C$ we then have the decomposition $A' \simeq \prod_{\sigma} A_{\sigma}(\C)$, where the product runs over all the embeddings $\sigma : K \hookrightarrow \C$ and $A_{\sigma}$ is the abelian variety obtained by action of $\sigma$ on $A$.  Then, $A'$ is of dimension $g' = g [K:\rat]$. Furthermore 
\begin{equation}\label{A'-A}
L(A'/\rat, s) = L(A/K,s),  \esp  \F_{A'/\rat} = N_{K/\rat}(\F_{A/K})\cdot D_{K}^{2g}, \textrm{ and } \Lambda(A'/\rat, s) = \Lambda(A/K,s).
\end{equation} 

The Hasse-Weil bound gives $P_{A', p}(T) = \prod_{i= 1}^{\rho}(1-\alpha_{i}T)$, where $\rho = \deg (P_{A', p}) \leq 2g'$ and $|\alpha_{i}| \leq \sqrt{p}$. Then, if we write $s = \sigma + i \tau$, with $\sigma = \Re(s) > \frac{3}{2}$, the local factor of the Eulerian product of the $L$-series satisfies
$$|P_{A',p}(p)^{-s}|^{-1} \leq (1-p^{\frac{1}{2} - \sigma})^{-2g'},$$
hence 
$$|L(A'/\rat, s)| \leq \zeta(\sigma-\frac{1}{2})^{2g'}.$$
Let $\sigma = \frac{3}{2} + \epsilon$, with $\epsilon > 0$. Let's denote $\F' = \F_{A'/\rat}$. Then
$$
|\Lambda(A'/\rat, s)| = |\Lambda(A'/\rat, \frac{3}{2} + \epsilon + i \tau)| \leq \F'^{\frac{3}{4} + \frac{\epsilon}{2}}\cdot (2 \pi)^{-g'(\frac{3}{2} + \epsilon)} \cdot \Gamma(\frac{3}{2} + \epsilon)^{g'} \cdot |\zeta(1+\epsilon)|^{2g'}.
$$
Since $|\zeta(1+\epsilon)| \leq (1+\frac{1}{\epsilon})$, for $\epsilon >0$, it follows that 
\begin{equation}\label{borne-lambda}
 |\Lambda(A'/\rat, \frac{3}{2} + \epsilon + i \tau)| \leq M(\epsilon), 
\end{equation}
with $M(\epsilon) = \F'^{\frac{3}{4} + \frac{\epsilon}{2}}\cdot (2 \pi)^{-g'(\frac{3}{2} + \epsilon)} \cdot \Gamma(\frac{3}{2} + \epsilon)^{g'} \cdot (1+\frac{1}{\epsilon})^{2g'}$.

Using the functional equation, that is, Conjecture \ref{funct-eq}, the same bound (\ref{borne-lambda}) is valid for $|\Lambda(A'/\rat, \frac{1}{2} - \epsilon - i \tau)|$. 

We now apply Lemma \ref{phragmen-lindeloef-lemma-strip}  to the function $\Lambda (A'/\rat, s)$, which is regular into the strip between the two parallel lines $\sigma = 3/2 + \epsilon$ and $\sigma = 1/2 - \epsilon$ and satisfies the bound (\ref{borne-lambda}) on these lines. Since we supposed that $\Lambda(A'/\rat, s) = \Lambda(A/K, s)$ satisfies Hypothesis \ref{hyp-phragmen-lindeloef}, a suitable choice for $\epsilon$ (made at the end of our proof) makes that condition (\ref{condition-phragmen-lindeloef-strip}) is satisfied. We then conclude that  the bound (\ref{borne-lambda}) is still valid throughout the strip, that is, for $s$ with real part $\sigma$ satisfying $\frac{1}{2}- \epsilon \leq \sigma \leq \frac{3}{2} + \epsilon$.

Applying the Cauchy inequality in the disc $\mathcal{D}(1, \frac{1}{2} + \epsilon)$, we obtain
$$L^{*}(A'/\rat, 1) = %
\frac{(2\pi)^{g'}}{\sqrt{\F'}} \frac{\Lambda^{(r)}(A'/\rat, 1)}{r!}
\leq \frac{(2\pi)^{g'}}{\sqrt{\F'}} \frac{1}{(\frac{1}{2} + \epsilon)^{r}} \max_{s \in \mathcal{D}(1, \frac{1}{2} + \epsilon)} \Lambda(A'/\rat, s).$$
The upper bound (\ref{borne-lambda}) gives
$$
L^{*}(A'/\rat, 1)  \leq \frac{1}{(\frac{1}{2} + \epsilon)^{r}} \cdot (2\pi)^{-g'\left(\frac{1}{2}+\epsilon \right)} \cdot \F'^{\frac{1}{4}+ \frac{\epsilon}{2}} \cdot \Gamma\left(\frac{3}{2} +\epsilon \right)^{g'} \cdot \left(1+\frac{1}{\epsilon}\right)^{2g'}.
$$
 To prove (\ref{bound-coeff-rank}), we choose $\epsilon = \frac{1}{2}$ and obtain 
$$
L^{*}(A'/\rat, 1)  \leq \left(\frac{9}{2\pi}\right)^{g'}\cdot \sqrt{\F_{A'/\rat}}.
$$
To prove (\ref{bound-coeff-cond}), we take $\epsilon = \frac{2}{\log \F'}$. (Remark that, since $A'$ is defined over $\rat$, $F_{A'/\rat} > 3$ and $\log F_{A'/\rat}  \ne 0$.) Thus $(\frac{1}{2} + \epsilon)^{-r} \leq 2^{r}$ and $\F'^{\frac{1}{4}+ \frac{\epsilon}{2}} = e\cdot \F'^{\frac{1}{4}}$.  
Since the abelian variety $A'$ is defined over $\rat$ and satisfies Conjecture \ref{funct-eq}, we have \cite[Section 3, Proposition]{mestre1986Compo}
\begin{equation}\label{N>10^{g}}
\F'= \F_{A'/\rat} > 10^{g'}. 
\end{equation}
Then, $1+\frac{1}{\epsilon} \leq \log \F'$ and $\frac{1}{2}+\epsilon \in [\frac{1}{2}, \frac{3}{2}]$, and therefore  
$ \Gamma\left(\frac{3}{2} +\epsilon \right) 
= (\frac{1}{2}+\epsilon)\Gamma\left(\frac{1}{2} +\epsilon \right)
\leq \frac{3}{2}\sqrt{\pi} < 3$.
Moreover $ 3^{g'}\cdot (2\pi)^{-g'(1/2 +\epsilon)} = \left(\frac{3}{(2\pi)^{1/2 + \epsilon}}\right)^{g'} \leq (\frac{6}{5})^{g'}$. This gives
$$
L^{*}(A'/\rat, 1)  \leq 2^{r} \cdot e \cdot (6/5)^{g'}\cdot \F_{A'/\rat}^{\frac{1}{4}}\cdot (\log \F_{A'/\rat})^{2g'}. \footnote{For avoiding the use of Conjecture \ref{funct-eq} a second time, instead of (\ref{N>10^{g}}) we can use $\F' > 3$ and obtain $L^{*}(A'/\rat, 1)  \leq e \cdot 2^{r} \cdot 5^{g'}\cdot \F_{A'/\rat}^{\frac{1}{4}}\cdot (1 + \frac{1}{2}\log \F_{A'/\rat})^{2g'}$.}
$$
We conclude both cases by applying (\ref{A'-A}). It remains to be shown 
that the hypothesis (\ref{condition-phragmen-lindeloef-strip}) on the growth of the $\Lambda$-series is satisfied by both choices of $\epsilon$. In fact, for $\epsilon = \frac{1}{2}$, $\alpha(\epsilon) = \frac{\pi}{2}$. For  $\epsilon = \frac{2}{\log \F'}$, using that $\F' > 2$, we obtain $\alpha(\epsilon) = \pi \times (\frac{4}{\log 2} + 1)^{-1}\leq \frac{\pi}{2}$. 
\hfill $\Box$

\medskip


\subsection{Bound for the local factor}\label{section-local-periods}

Is in this section where we assume  that $A$ is principally polarised. For Lemma  \ref{c_v-real}, this hypothesis is technical. It can be removed from Theorem \ref{autissier-matrix-lemma}, as in the version given in \cite[Th\'eor\`eme 1.1]{gaudron-remond-periodes2014}.   Nevertheless, the exposition  is simplified here, and we also use it in order to make evident the $\Im \tau_v$  in the matrix lemma (Theorem \ref{autissier-matrix-lemma} and Lemma \ref{rho-im-tau} and \ref{matrix-lemma}). 

We will bound the local periods $c_v$. For every non-archimedean place $v$, the  numbers $c_v$ are non-zero integers and can be bounded from below by 1.

As for the archimedean local periods, in order to relate them to the Faltings' height, we need some preliminaries.

For $v$ complex, the local period $c_{v}$  is almost the norm $||\omega||_{v}$ of $\omega$ (up to the normalisation of the Hermitian metric, see Lemma \ref{c_v-complex}), while for $v$ real, it is a little bit more delicate to link the local period $c_{v}$ with the norm $||\omega||_{v}$ (see Lemma \ref{c_v-real}).

We fix an archimedean place $v$ of $K$. Recall that  $(\gamma_{1,v}, \ldots, \gamma_{2g,v})$ is a basis of the integral homology $H = H_1(A(\overline{K_{v}}), \Z)$ of $A$, chosen so that $\gamma_{1,v}, \ldots, \gamma_{g,v}$ generates the part of $H$ fixed by complex conjugation. Let 
\begin{equation}\label{Omegas}\Omega_{1,v} = \left(\int_{\gamma_{i,v}} \omega_{j}\right)_{1\leq i \leq g} \esp \textrm{and} \esp   \Omega_{2,v} = \left(\int_{\gamma_{i,v}} \omega_{j} \right)_{g+1\leq i \leq 2g}
\end{equation}
be the periods matrixes associated to $\gamma_{1,v}, \ldots,\gamma_{2g,v}$, where $j$ runs over $\{1, \ldots, g\}$. Moreover, since $A$ is principally polarised, 
we can choose $\gamma_{1,v}, \ldots\gamma_{2g,v}$ such that 
$$\tau_{v} = \Omega_{1,v}^{-1}\Omega_{2,v}$$
is a symmetric matrix in the Siegel space, that is, $\Im \tau_v$ is definite positive. Let $\Lambda_{v} = \Omega_{1,v} \Z ^{g}+ \Omega_{2,v} \Z^{g}$ be the associated lattice. %
Choose an isomorphism over $\C$ 
\begin{displaymath}
\begin{array}{rcl}
\varphi : \C^{g}/\Lambda_{v} &  \rightarrow & A(\overline{K_{v}}) %
\end{array}
\end{displaymath}
such that the inverse function of $\varphi$ maps the invariant differential $\eta = \omega_{1}\wedge \ldots \wedge \omega_{g}$ to $dz$: 
$$\varphi^{*}(\eta) =  dz.$$

Let $\Gamma_{v} = \Z ^{g}+ \tau_v  \Z^{g} = \Omega_{1,v}^{-1} \Lambda_{v}$ and choose also an isomorphism over $\C$ 
\begin{displaymath}
\begin{array}{rcl}
\psi : \C^{g}/\Gamma_{v} &  \rightarrow & A(\overline{K_{v}}) %
\end{array}
\end{displaymath}
such that
$$\psi^{*}(\eta) =  \det \Omega_{1,v} dz.$$


 We deduce the next result from Rémond's work \cite{these.remond}. 
 
\begin{lemma}\label{c_v-complex} 
For a complex place $v$, we have 
\begin{equation}
c_v = 2^g ||\eta||_{v}.
\end{equation}
\end{lemma}

\noindent {\it Proof.}
We use the module measure of the differential form $\eta$, (\ref{mod-measure-archimedean}) and (\ref{mod-measure-wedge-complex}):

$$c_v := \left | \det\left( \int_{\gamma_i} \omega_i \right)_{1 \leq i,j \leq 2g} \right|= \int_{A(K_v)} \mod(v) = \int_{A(\C)} |\eta \wedge \overline{\eta}| =: 2^g  ||\eta||_{v} .$$
\hfill$\Box$\medskip

\begin{lemma} \label{c_v-real} 
We suppose  that $A$ carries a principal polarisation.  For a real  place  $v$, we have 
$$c_{v} =  \frac{|A_{v}(\R) : A_{v}(\R)^{0}|}{\sqrt{\det \Im(\tau_{v})}} ||\eta||_{v}.$$
\end{lemma}

\noindent {\it Proof.}
By definition of $c_v$ and $\Omega_{1,v}$ we have
$$c_{v} := |A_{v}(\R) : A_{v}(\R)^{0}| \cdot \left|\det \left(\int_{\gamma_{i,v}} \omega_{j}\right)_{1\leq i,j\leq g} \right| =: |A_{v}(\R) : A_{v}(\R)^{0}| \cdot  |\det \Omega_{1,v}|.$$
Then, using the definition of the metric for $v$ real and the inverse map of $\psi$ we compute
 $$|| \eta ||_{v}^2 = |\eta|_{v}^2 :=  \frac{1}{2^g} \int _{A(\overline{K_{v}})} \left | \eta \wedge \overline{\eta}   \right | =  \frac{1}{2^g}  \int_{\C^{g}/\Gamma_{v}}  |\det \Omega_{1,v}|^{2} \left | dz \wedge \overline{dz}  \right | $$
$$ =  \frac{1}{2^g}   |\det \Omega_{1,v}|^{2} \int_{\C^{g}/(\Z^{g} + \tau_{v} \Z^{g})}  2^g  \left |  dx \wedge dy  \right | =  |\det \Omega_{1,v}|^{2} \det \Im(\tau_{v}).$$
For the last equality we use that $\int_{\C^{g}/(\Z^{g} + \tau_{v} \Z^{g})}  \left | dx \wedge dy \right | $ is the area of a fundamental domain for $\C^{g}/(\Z^{g} + \tau_{v} \Z^{g}$), which is $\det \Im (\tau_v)$.
Then
$$ || \eta ||_{v} = \frac{ \sqrt{\det\Im (\tau_{v})}}{|A_{v}(\R) : A_{v}(\R)^{0}|} c_v.$$
\hfill$\Box$\medskip

We then deduce the following result for $c_{\infty}(A/K)$.

\begin{lemma} \label{local-falt} We assume  that $A$ carries a principal polarization. 
The archimedean local factor defined in (\ref{archimedean-local-factor}) satisfies 
$$c_{\infty}(A/K) = 2^{g t} \prod_{v {\textrm{real}}}\frac{|A_{v}(\R) : A_{v}(\R)^{0}|}{\sqrt{\det \Im(\tau_{v})}} \cdot  e^{-[K:\rat] h_{Falt}(A/K)},$$
where $t$ is the number of complex places of $K$. 
\end{lemma}

\medskip

\noindent {\it Proof.} 
 Recall that $H^{0}(\A, \DiffA) = \eta \cdot \va_{\eta}$. Recall also that,  by the product formula, $\sum_{v \in M_{K}}\log ||\eta||_{v} = \sum_{v \in M_{K}}\log ||k \eta||_{v} $, for all $k$ in $K$. Then to compute the degree of $\eDiffA$, we choose the invariant differential $\eta$:

$$[K:\rat] h_{Falt}(A/K)=\deg_{Ar}(\eDiffA, ||.||) = \log |\eta \va_{\eta}/\eta \ok| - \sum _{v |\infty} \log ||\eta||_{v}.$$

On the one hand, $|\eta \va_{\eta}/\eta \ok| = | \va_{\eta}/ \ok| = |\ok/\va_{\eta}^{-1}| = N_{K/\rat}(\va_{\eta}^{-1})$.

On the other hand, from  Lemma \ref{c_v-complex}  and  Lemma \ref{c_v-real}, we deduce that the product of the archimedean local periods satisfies the following equality
$$\prod_{v|\infty} ||\eta||_{v}=  \prod_{v \,\textrm{real}}\frac{ \sqrt{\det \Im(\tau_{v}}) }{|A_{v}(\R) : A_{v}(\R)^{0}|} \cdot \prod_{v \,\textrm{complex}}\frac{ 1}{2^g} \cdot \prod_{v|\infty} c_v,$$
and then
\begin{equation} \label{Ab-var-height-local-Imtau}
[K:\rat] h_{Falt}(A/K)= \log N_{K/\rat}(\va_{\eta}^{-1}) -   \sum_{v|\infty} \log c_v -   \sum_{v \,\textrm{real}}  \log \frac{ \sqrt{\det \Im(\tau_{v}}) }{|A_{v}(\R) : A_{v}(\R)^{0}|} + \sum_{v \,\textrm{complex}}  \log (2^g).
\end{equation}
Since, by definition (\ref{archimedean-local-factor}),  $c_{\infty}(A/K) = N_{K/\rat}(\va) \cdot \prod_{v \in M_K^{\infty}}c_v$,   we can conclude.
\hfill$\Box$\medskip

\medskip

Remark that it is the Faltings' height on $K$, instead of the {\it stable} one, which appears naturally in the proof.

It is worth noting  that Rémond \cite[Lemme II.3.1]{these.remond} proves an analog result to our Lemma \ref{c_v-real}, without the assumption that the polarisation is principal.  

Nevertheless, in order to apply the matrix lemma and relate the local periods with the Faltings' height, we will assume later that $A$ carries a principal polarisation.


To state his result, let's fix an archimedean place $v$ of $K$.  Denote $\Lambda  = H_1(A_v(\overline{K_{v}}), \Z)$ the integral homology  of $A_v$, $\Lambda^+$ its sub-module fixed by complex conjugation, and $\Lambda^-$ the biggest sub-module of $\Lambda$ in which the complex conjugation induces $-id$. 
Let's fix $(\gamma_{1,v}, \ldots, \gamma_{g,v})$ a basis for  $\Lambda^+$, and $(\gamma_{g+1,v}, \ldots, \gamma_{2g,v})$ a basis for  $\Lambda^-$. Let $M_v$ be a matrix of $GL_g(\C)$ such that
$$(\gamma_{g+1,v}, \ldots, \gamma_{2g,v}) = M_v (\gamma_{1,v}, \ldots, \gamma_{g,v}).$$
Since $M_v$ depends on the choice of the basis up to a multiplication by an element of $GL_g(\enteros)$,  its determinant is well defined up to a sign. Let's choose the basis such that $\det (\Im M_v)>0$.

Then the  transcription of  \cite[Lemme II.3.1 p.18]{these.remond} with our notations and normalisations reads as 
\begin{equation}\label{c_v-real-Gael-M_v} 
c_{v} = 2^{\frac{g}{2}}  \frac{|A_{v}(\R) : A_{v}(\R)^{0}|^{\frac{1}{2}}}{\sqrt{\det \Im(M_{v})}} ||\eta||_{v},
\end{equation}
for a real place $v$, and his Corollaire II.3.1 gives
\begin{equation}\label{local-falt-Gael} 
c_{\infty}(A/K) =  \prod_{v {\textrm{real}}} \frac{|A_{v}(\R) : A_{v}(\R)^{0}|^{\frac{1}{2}}} {\sqrt{\det \Im(M_{v})}} \cdot 2^{\frac{g[K:\rat]}{2}} \cdot e^{-[K:\rat] h_{Falt}(A/K)}.
\end{equation}

In order to compare our results with Rémond's ones, recall that, we have chosen $\tau_{v} = \Omega_{1,v}^{-1}\Omega_{2,v}$, where $\Omega_{1,v}$ and $\Omega_{2,v}$ are the periods matrixes associated to $\gamma_{1,v}, \ldots,\gamma_{2g,v}$ defined by (\ref{Omegas}),  being a symmetric matrix in a  fundamental domain.

We then have an induced isomorphism $A(\overline{K_{v}})  \simeq \C^g/\Gamma_v$ with $\Gamma_v = \enteros^g + \tau_v \enteros^g$ and, with the above notation, $\Lambda \simeq \Gamma_v = \enteros^g + \tau_v \enteros^g$. Then $\Lambda^+ \simeq \enteros^g$ and $\Lambda^- \simeq M_v \enteros^g$. Moreover $\Lambda^+ \oplus \Lambda^- = \enteros^g + M_v \enteros^g \subset \Lambda = \Gamma_v = \enteros^g + \tau_v \enteros^g$. 

On the one hand, the index satisfies the following equalities.
$$[\Lambda : \Lambda^+ \oplus \Lambda^-] = \frac{\Vol(\Lambda)}{\Vol (\Lambda^+ \oplus \Lambda^-)} = \frac{\Covol(\Lambda^+ \oplus \Lambda^-)}{\Covol(\Lambda)} = \frac{\det \Im(M_v)}{\det \Im(\tau_v)}.$$

On the other hand, we have (see the proof of \cite[Lemme II.3.1]{these.remond}),
$$[\Lambda : \Lambda^+ \oplus \Lambda^-] = \frac{2^g}{|A_{v}(\R) : A_{v}(\R)^{0}|}.$$

This gives 
$$\det \Im \tau_v = \frac{|A_{v}(\R) : A_{v}(\R)^{0}|}{2^g} \det \Im M_v.$$
(See \cite[(A.3) p. 79]{these.remond} for explicit  $ \tau_v$ in the one-dimensional case.) We then deduce from (\ref{c_v-real-Gael-M_v}) the same result for $c_v$, when $v$ is real, as our Lemma \ref{c_v-real}, and from (\ref{local-falt-Gael}) the same result for $c_{\infty}(A/K)$ as our Lemma \ref{local-falt}.

Observe that the term $2^g$, and respectively $2^{\frac{g}{2}} $, in Lemma \ref{c_v-complex}, and, respectively in (\ref{c_v-real-Gael-M_v}), come from our choice of the Hermitian metric (\ref{norme2}), in which we choose to use Faltings' normalisation $\frac{1}{2^g}$.  


However, the $c_{\infty}(A/K)$, and, by Lemma \ref{local-falt}, the period matrices $\Im \tau_v$, appear in the BSD-formula. Since we would like to bound it by more tractable objects associated to out variety, we will make use of a {\it matrix lemma}. 

We call a matrix lemma an upper bound for the period matrix in terms of the height of the abelian variety. Such a relation was first introduced by Masser \cite[p. 115]{masser.LNM1290}. See also \cite[Lemma 8.6 p.440]{masser-wuestholz-periods1993}. A new approach in terms of the Faltings' height was introduced by  Bost \cite{zbMATH00908680}, \cite{zbMATH01003129}. Further  and effectives versions are due to  Graftieaux \cite{graftieaux2001},  David-Philippon \cite[Lemma 6.7]{sinnou-patrice.helvet}, and  Gaudron \cite{gaudron-ENS-2006}. 
We use here  Autissier's result \cite[Corollaire 1.4]{autissier-matrixlemma}, weakened because stated with the Faltings' height over $K$ instead of the stable one.

Denote $(A, \mathcal{L})$ the abelian variety carrying the principal polarisation $\mathcal L$, $H$ the Riemann form associated to $\mathcal{L}$, $t_{A}$ the tangent space of $A$ at the origin, and $\Omega_{A}$ its period lattice.  The form $H$ is an hermitian form on $t_A$ definite positive satisfying $\Im H(\Omega_A, \Omega_A) \subset \enteros$. The polarisation $\mathcal{L}$ endow the tangent space $t_A$  with an hermitian norm. We define, for $z \in t_A$,  $||z||_{\mathcal{L}}^2 := H(z,z)$. For an archimedean place $v$, we have the same objects related to $A_v$, which we denote with a sub-index $v$. That is, we denote $||.||_{\mathcal{L}_v}$ the hermitian norm induced  by the polarisation $\mathcal{L}_v$ into the tangent space  $t_{A_v}$.  Let $\rho(A_v, \mathcal{L}_v)$ be the minimum value of $||\omega||_{\mathcal{L}_v}$, for all non-zero $\omega$ in $\Omega_{A_v}$.

\begin{thm}[Autissier]\label{autissier-matrix-lemma} Suppose that the abelian variety $(A, \mathcal{L})$ carries a principal polarisation. Then
$$\frac{1}{[K:\rat] }\sum_{v \in M_K}(\rho(A_v, \mathcal{L}_v))^{-2} \leq 3 h_{Stab}(A) + 6g.
$$
\end{thm}

\begin{lemma}\label{rho-im-tau}
For $\mathcal L$ a principal polarisation and $v$ an archimedean place, we have
$$\rho(A_v, \mathcal{L}_v))^{-2} \geq \frac{1}{g} (\det \Im \tau_v )^{1/g}.$$
\end{lemma}

\medskip

\noindent {\it Proof.} We choose to use the parametrization $A(\overline{K_{v}})  \simeq \C^g/\Gamma_v$, with $\Gamma_v = \enteros^g + \tau_v \enteros^g$. Then, for $z \in \C^g$, we denote $ ||z||_{\mathcal{L}_v}^2 =  ||z||_{\Gamma_v}^2 = {H(z,z)} = {}^t \! z (\Im \tau_v )^{-1} \bar{z}$. 

We apply Minkowski's theorem to the lattice $\Gamma_1 = \enteros^g \subset \Gamma_v$, that is,   $\lambda_1(\Gamma_1)^g \Vol(\Gamma_1) \leq 2^g$, where $\lambda_1(\Gamma_1)$ is the first minimum of $\Gamma_1$ for the norm $ ||z||_{\Gamma_v}$. 
We then have
$$\lambda_1(\Gamma_1)^g \leq  2^g \frac{\Covol(\Gamma_1)}{\Vol (B)},$$
where $B$ is the unit ball. On one hand, for the chosen norm we have $\Covol(\Gamma_1) =  (\det \Im \tau_v)^{-1/2}$. On the other hand, the unit ball $B$ contains 
$[-\frac{1}{\sqrt{g}}, \frac{1}{\sqrt{g}}]^g$,  and then $\Vol(B)  \geq \frac{2^g}{g^{g/2}}$.  
Hence
$$\lambda_1(\Gamma_1)^g \leq \frac{g^{g/2}}{\sqrt{\det \Im \tau_v}}.$$
Since $\rho(A_v, \mathcal{L}_v):= \lambda_1 (\Gamma_v) \leq \lambda_1(\Gamma_1)$, we then have
$$ \rho(A_v, \mathcal{L}_v) \leq \frac{g^{1/2}}{(\det \Im \tau_v)^{\frac{1}{2g}}},$$
and we can conclude.
\footnote{Remark that we could relate this with the parametrization we used before, that is $A(\overline{K_{v}})  \simeq \C^g/\Lambda_v$, with $\Lambda_v = \Omega_{1,v} \enteros^g + \Omega_{2,v} \enteros^g$. Indeed, we could define   $||z||_{\Lambda_v}^2  = H'(z,z) = {}^t \! z  {}^t  \Omega_1^{-1} (\Im \tau_v)^{-1}  \overline{\Omega_{1,v}^{-1} } \overline{z}$, 
since $\Im{H'(\Lambda_v, \Lambda_v)} \subset \enteros$.
 Denote $ \lambda_1' (\Omega_{1,v} \enteros^g)$ the first minimum of the lattice $\Omega_{1,v} \enteros^g$ for the norm $||z||_{\Lambda_v}^2$. 
 We then have $ \lambda_{1}(\enteros^g) = \lambda_1' (\Omega_{1,v} \enteros^g)$.}
\hfill$\Box$\medskip

We then deduce from Theorem \ref{autissier-matrix-lemma} the next matrix lemma, involving $\tau_v$. 

\begin{lemma}\label{matrix-lemma}
Suppose that the abelian variety $(A, \mathcal{L})$ carries a principal polarisation. The sum of the determinants of the matrices $\Im(\tau_{v})$ satisfies
$$\frac{1}{g [K:\rat] }\sum_{v \in M_K}(\det \Im(\tau_{v}))^{1/g} \leq 3 h_{Stab}(A) + 6 g.
$$
\end{lemma}

{Nevertheless, in the case of an elliptic curve, we could work  directly  with the Faltings' height. 
  In fact, for $A$ an elliptic curve, we have (\cite[Prop. 1.1 of Chap. X]{cornell-silverman})
\begin{equation}\label{Ellip-c-height-local-Imtau}
12[K:\rat] h_{Falt}(A/K) = \log N_{K/\rat} \Delta_{A/K} - \sum_{v|\infty} n_v \log |\Delta(\tau_v)| - \sum_{v | \infty} 6 n_v \log \Im (\tau_v),
\end{equation}
where $\Delta_{A/K}$ is the minimal discriminant of the elliptic curve, and $\Delta (\tau)$ is the modular form  $(2 \pi ) ^{12} q_{\tau} \prod_{n = 1}^{\infty} (1 - q_{\tau}^n)^{24}$, where $q_{\tau} = e ^{2 \pi i \tau}$. 

It is worth noting that this is an analogous formula to (\ref{Ab-var-height-local-Imtau}). On the right hand term, the first term involves an ideal of $K$, that is $\va_{\eta}$, and, resp. $\Delta_{A/K}$. The second term of (\ref{Ab-var-height-local-Imtau}) involves the uniformization $\varphi : \C^{g}/\Lambda_{v}   \rightarrow  A(\overline{K_{v}})$. In fact, by definition, for the real archimedean places, the $c_v$ are related to the matrix  $\Omega_{1,v}$ and $\Lambda_{v} = \Omega_{1,v} \Z ^{g}+ \Omega_{2,v} \Z^{g} = \Omega_{1,v} (\Z ^{g}+ \tau_{v} \Z^{g} )$. In formula (\ref{Ellip-c-height-local-Imtau}), $\Delta(\tau_v)$ is the discriminant of the equation of the curve defined by the parametrization $ \C/(\Z + \tau_{v}  \Z ) \rightarrow  A(\overline{K_{v}})$. Both third terms involve $\Im (\tau_{v})$. Observe that, while in (\ref{Ellip-c-height-local-Imtau}) all the archimedean places are involved in this third term, in (\ref{Ab-var-height-local-Imtau}), only the real ones appear, and some parts of the $\tau_{v}$ are hidden in the second term. This is because, in the second term of (\ref{Ab-var-height-local-Imtau}), for the real places, the $c_{v}$ are related to $\Omega_{1,v}$, and the parametrization is the one related to the lattice  $\Omega_{1,v} (\Z ^{g}+ \tau_{v} \Z^{g} )$, while in (\ref{Ellip-c-height-local-Imtau}), the lattice is $\Z + \tau_{v}  \Z $.  (The fourth term is there because the normalisation of the norm that we use to define Faltings' height.)
However, the proof of both formulas are close, and they involved the area of the fundamental domain of the lattice. 

Moreover,  if we work out the third term (the one with the $\Im (\tau_{v})$), e.g. with a matrix lemma, we could obtain a relation between the height of the variety $h_{Falt}(A/K)$ and the local numbers $c_{\infty}(A/K)$. That is what we will do in Lemma \ref{lemme-local-falt-ineg}  in general dimension. 
In fact, we deduce from Lemma \ref{local-falt} a lower bound for $c_{\infty}(A/K)$, in terms of the Faltings' height,  the degree $[K:\rat]$ and the dimension $g$. It is this lower bound,  which we will use to prove the main results.

Conversely,  we point out that, if we work out the second term (which is $|\Delta(\tau_v)|$ for (\ref{Ellip-c-height-local-Imtau})), we could obtain a relation between $h_{Falt}(A/K)$ and the $\Im (\tau_{v})$, that is, a {matrix lemma}.
That is what we will  do in Lemma \ref{matrix-lemma-g=1-Andrea}, where we recover the same kind of result as in Lemma \ref{matrix-lemma},  in the one-dimensional case.

\begin{lemma}\label{matrix-lemma-g=1-Andrea}
If $g=1$, the following bound also holds
$$\frac{1}{[K:\rat]}\sum_{v \in M_K}\Im(\tau_{v}) \leq 3 h_{Falt}(A/K) + \frac{11}{2}. 
$$
\end{lemma}

\medskip

\noindent {\it Proof.}
We use the formula (\ref{Ellip-c-height-local-Imtau}). 

We estimate then each of the three terms on the right hand side. First, $\log N_{K/\rat} \Delta_{A/K} \geq 0$. Second, $\log \Im (\tau_v)  \leq \frac{1}{e} \Im (\tau_v)$, hence $- \sum_{v | \infty} 6 n_v \log \Im (\tau_v) \geq -\frac{6}{e}  \sum_{v | \infty}  n_v \Im (\tau_v)$. As for the second term, we use the estimate (see exercise on page 256 of {\it loc. cit.}, where we corrected the missprint for the modular form),
$$\log |\Delta(\tau_v)| = \log \left | (2 \pi)^{12} q_{\tau} \right | + A_{\tau}, {\rm whith} \left | A_{\tau} \right | \leq \frac{1}{9}.$$ 
Thus 
$$\log |\Delta(\tau_v)| \leq -2 \pi \Im (\tau_v) + \log \left( (2\pi)^{12}e^{1/9} \right) . $$
 We obtain
$$\frac{1}{[K:\rat] } \sum_{v|\infty}  \Im (\tau_v)  \leq \frac{1}{[K:\rat] }  \sum_{v|\infty} n_v  \Im (\tau_v)  \leq \frac{12 }{2\pi - 6/e} h_{Falt}(A/K) + \frac{ \log  \left( (2\pi)^{12}e^{1/9} \right) }{2\pi - 6/e},$$
what gives the announced bound.
\hfill$\Box$\medskip
\medskip

Observe that Proposition 3.2 of \cite{gaudron-remond-periodes2014} (which follows from Autissier's matrix lemma) with Remarque 3.3 of {\it loc. cit.} gives bounds in the case where $A$ is an elliptic curve. 
They use Deligne's normalisation for the stable Faltings' height, which we denote by $h_D(A)$, and with our notation this reads as $h_D(A) = h_{Stab}(A) + g/2 \log \pi$. Using the fact that $h_{Stab}(A) \leq h_{Falt}(A/K) $, we deduce the following bounds.

\begin{lemma}\label{matrix-lemma-g=1}
If $g=1$, then 
$$[K:\rat]^{-1} \sum_{v|\infty}  \Im (\tau_v) \leq 6,45 \max\{h_{Falt}(A/K) + \frac{\log \pi}{2},1  \}$$
and also 
$$[K:\rat]^{-1} \sum_{v|\infty}  \Im (\tau_v) \leq  1,92\max\{h_{Falt}(A/K) + \frac{\log \pi}{2},1000  \}.$$

\end{lemma}

Observe that for $h_{Falt}(A/K) \in [\frac{1}{2}, 638]$, the bound of Lemma \ref{matrix-lemma-g=1-Andrea} is sharper than the bounds of Lemma \ref{matrix-lemma-g=1}.
In contrast, for $h > 638$, the second bound of Lemma  \ref{matrix-lemma-g=1} is of a better quality. In particular for $h \in [638, 1000- \frac{\log \pi}{2}]$, where it equals 1920, which is independent of the height. Moreover, for $h > 1000- \frac{\log \pi}{2}$, the bound $3, 02 h$ is valid too. 
However, all of these bounds have the shape $$\frac{1}{[K:\rat]}\sum_{v \in M_K}\Im(\tau_{v}) \leq m_1 h_{Falt}(A/K) + m_2,
$$
with $m_1 \leq 6,45$ and $m_2 \leq 1920$. In what follows, to simplify the exposition and because it doesn't matter for our purpose, we will use $3h + 6$, which always holds. 
\medskip


\begin{lemma}\label{lemme-local-falt-ineg}
The archimedean local factor satisfies the following inequality.

\begin{equation}\label{bound-local-falt-ineg}
 c_{\infty}(A/K) ^{-1} \leq (3 g [K:\rat]  h_{Falt}(A/K) + 6 g^2 [K:\rat]) ^{\frac{g[K:\rat]}{2}}\cdot e^{[K:\rat] h_{Falt}(A/K)}. 
\end{equation}
The weaker, but eventually more useful bound also holds.
\begin{equation}\label{local-falt-ineg-weaker}
 c_{\infty}(A/K) ^{-1} \leq   (6g^2[K:\rat] \max\{1, h_{Falt}(A/K)\})^{\frac{g[K:\rat]}{2}}\cdot e^{[K:\rat] h_{Falt}(A/K)}. 
\end{equation}
\end{lemma}

\smallskip

\noindent {\it Proof.} 

Let's denote $d= [K:\rat]$ and $t$ the number of complex places of $K$.  Since $|A_{v}(\R) : A_{v}(\R)^{0}| \geq 1$,  we deduce from Lemma  \ref{local-falt} 
$$c_{\infty}(A/K)^{-1} \leq \left(\frac{1}{2}\right)^{gt} \prod_{v \in M_{K}^{\infty}{\textrm{real}}}\sqrt{\det \Im(\tau_{v})}\cdot e^{d h_{Falt}(A/K)}. 
$$
For the complex places, we use that  $\det \Im(\tau_{v}) \geq (\frac{\sqrt{3}}{2})^g$ (because $\tau_v$ is a matrix in the Siegel space). Thus 
$$ \prod_{v \in M_{K}^{\infty}{\textrm{real}}}\sqrt{\det \Im(\tau_{v})} \leq \left(\frac{2}{\sqrt{3}}\right)^{\frac{gt}{2}} \prod_{v \in M_{K}^{\infty}}\sqrt{\det \Im(\tau_{v})}.$$
 Using the arithmetic-geometric inequality, we obtain 
$$\prod_{v \in M_{K}^{\infty}} (\det \Im (\tau_v))^{1/g} \leq \frac{1}{d} \left(\sum_{v \in M_{K}^{\infty}} (\det \Im (\tau_v))^{1/g}\right)^{d}.$$
And then, 
$$c_{\infty}(A/K)^{-1} \leq \left(\frac{1}{2}\right)^{gt} \left(\frac{2}{\sqrt{3}}\right)^{\frac{gt}{2}}   \left(\sum_{v \in M_{K}^{\infty}} (\det \Im (\tau_v))^{1/g}\right)^{\frac{gd}{2}} \cdot e^{d h_{Falt}(A/K)}.$$
We conclude with Lemma  \ref{matrix-lemma}. 
The second inequality is just an easy deduction from the first one. 
\hfill$\Box$\medskip

\subsection{Bound for the cardinality of the torsion part}


In the one-dimensional case, using the results of Merel \cite{merel} and  Parent \cite{parent.torsion1999} we can obtain a uniform bound for the cardinality of the torsion part of the Mordell-Weil group. In fact, Merel's result tell us which prime numbers could divide $|E(K)_{\tors}|$ and Parent's result give us a bound for the powers of these primes, independent on the power.

\begin{lemma}\label{lemma.torsion.ellip}
For every integral number $d \geq 1$ there is a positive number $B({d})$ such that for every number field $K$ with $[K:\rat] \leq d$ and every elliptic curve $E$ defined over $K$ we have
\begin{equation}\label{merel'sbound}
|E(K)_{\tors}| \leq B({d}).
\end{equation}
On may take 
$B(d) =  (129. (5^d-1)(3d)^6)^{\frac{(1+3^{d/2})^8}{2\log (1+3^{d/2})}}$.
\end{lemma}

Notice that the dependence on $d$ of $B(d)$ is twice exponential and not only exponential as expected in \cite{parent.torsion1999}. 
For the convenience of the reader, we give the details  of the proof of Lemma \ref{lemma.torsion.ellip}. Before the proof, we state an analytic lemma, which will be used therein.

\begin{lemma}\label{lemma.analytic}
For $n \geq 1$, denote $p_1, p_2, \ldots, p_n$ the  $n$ first prime numbers. As usual, denote  $\theta(p_n)= \sum_{i=1}^n \log p_i$. For every  $n\geq 2$, one has 
$$ n \leq 4\, \frac{\theta(p_n)}{\log \theta(p_n)}.$$
\end{lemma}

\noindent {\it Proof.}
Remark (see, e.g.,  \cite[page 25]{ellison}) that for every $n\geq 1$, one has $p_n \geq n \log n$.
Furthermore, for $n \geq 2$, one has $\sum_{i=1}^n \log i \geq \int_1^n \log x dx = n \log n -n +1$ and $\sum_{i=2}^n \log (\log i) >\log \log 2$. From these remarks we deduce that, for $n\geq 2$,
$$\theta(p_n) = \log 2 + \sum_{i=2}^n \log p_i > \log 2 + \sum_{i=2}^n \log (i \log i) >\log 2 +  n \log n - n + 1 + \log \log 2.$$
Let $n\geq 4$. Then $\theta(p_n) > \frac{1}{2}\,n \log n \geq e$ and, since for $x \geq e$, the fonction $x \mapsto \frac{x}{\log x}$ is increasing,  then
$\frac{\theta(p_n)}{\log \theta(p_n)} \geq \frac{\frac{1}{2}\,n \log n}{\log(\frac{1}{2}\,n\log n)}$.
Moreover, $\frac{\log n + \log \log n - \log 2}{\log n}  = 1 + \frac{\log \log n}{\log n} - \frac{\log 2}{\log n} \leq 1 + \frac{1}{e}$. Thus 
$n \leq 2\,\left(1 + \frac{1}{e}\right) \frac{\theta(p_n)}{\log \theta(p_n)}$. 
We easily check that for $n =1, 2$ and $3$ one also has $n \leq 4 \frac{\theta(p_n)}{\log \theta(p_n)}$. 
\hfill $\Box$
\medskip

\noindent {\it Proof of Lemma \ref{lemma.torsion.ellip}.}  
Following a result of  Merel, if there is an element in $E(K)_{\tors}$ of order a prime number $p$, then $p \leq m(d)$. The theorem of \cite{merel} gives $m(d) = d^{3d^2}$; but this bound was improved by  Oesterl\'e (in an unpublished article) by $m(d) = (1+3^{d/2})^2$. We will use here Oesterl\'e's bound. 
Let us denote $p_1 < ... < p_m$ the first $m$ prime numbers, where $m$ satisfies $p_m \leq m(d)$ 
and $p_{m+1} > m(d)$. %
Since $m(d) \geq 4$, $m \geq4$, and $\theta(p_m) = \log (p_1... p_m) \geq \log (2 \times 3 \times 5 \times 7) \geq e$.
We also have $\theta(p_m)  \leq m \log m(d)$. %
Applying Lemma \ref{lemma.analytic} to $m$,  we deduce that 
$$ m \leq 
\frac{m(d)^4}{\log m(d)} = \frac{(1+3^{d/2})^8}{2\log (1+3^{d/2})}.$$
For $i \in \{1, \ldots , m\}$, there exist some $n_i \geq 0$, such that $|E(K)_{\tors}| \leq p_1^{n_1} ... p_m^{n_m}$. From  \cite[Theorem 1.2]{parent.torsion1999}, we know that, for every $p \in \{p_1, \ldots, p_m\}$ and every non-zero integer $n$,
$$p^n \leq c(d) = 129. (5^d-1)(3d)^6. $$
(In fact, Parent's result is even more precise; it gave  better bounds for $p^n$ depending if $p$ equals 2, 3 or not.) We conclude that 
$$|E(K)_{\tors}| \leq c(d)^m \leq (129. (5^d-1)(3d)^6)^{\frac{(1+3^{d/2})^8}{2\log (1+3^{d/2})}}.$$
\hfill$\Box$\medskip

In the general case, no such uniform bound on torsion of an abelian variety $A/K$, depending only on the dimension and the degree is known, but the following lemma suffice for our purpose.

\begin{lemma} \label{lemma-torsion}
 Let's  $\F = N_{K/\rat}(\mathcal{F}_{A/K})$ denote the norm of the conductor of $A/K$ and  $\G = \max\{2, \F\}$.  We have
$$|A(K)_{\tors}| \cdot |\Ad(K)_{\tors}| \leq \frac{5}{\log 2} \cdot (\log \G)^{4g[K:\rat]}.
$$
\end{lemma}

\noindent {\it Proof.}
As usual, let us denote $\omega(N)$ the number of prime numbers dividing $N$ and $\pi(X)$ the number of prime numbers $\leq X$. %
By \cite{dusart1999}, for $X \geq 17$, $\pi(X) \geq \frac{X}{\log X -1}$. And by \cite{robin1983}, for $N \geq 3$, we have  $\omega(N) \leq 1,3841 \frac{\log N}{\log \log N}$.

Set $Y = \log \F$ and  $C= 5/(\log 2)$. If $\F \geq \exp (\frac{17 \log 2}{5})$, which is bigger than 10, then $\F \geq 3$ and $CY = \frac{5}{\log 2} \log \F \geq 17$, and we can apply the two previous results, for $N = \F$ and $X= CY$. (In particular, $\log \log \F \ne 0$.) We then have

$$\pi(C \log \F) -\omega(\F)  = \pi(C Y) -\omega(e^Y) \geq \frac{CY}{\log (CY) -1} -1,3841 \frac{Y}{\log Y}$$
$$ \geq \frac{Y (C \log Y - 1,3841 \log Y + 1,3841 (1-\log C))}{(\log Y) (\log Y + \log C -1)} $$
$$\geq \frac{Y (C \log Y - 1,3841 \log Y)}{(\log Y) (\log Y)}  =  \frac{Y}{\log Y} (C - 1,3841) =   \frac{\log \F}{\log \log \F} \left(\frac{5}{\log 2} - 1,3841\right) \geq 2.$$

Moreover, we could verify (e.g. with GP/Pari), that there is no $\F \in [2, 11]$, for which the inverse inequality holds, that is, for all $\F \geq 2, \pi(\frac{5}{\log 2} \log \F) -\omega(\F) \geq  2$.

We can then take  two distinct primes numbers, $p$ and $q$, coprime with $\G= \max\{2, \F\}$ and $\leq \frac{5}{\log 2} \log \G$. (This is clear for if $\G = \F$. Otherwise $\F < 2$ and  we have $\{p, q \} \subset \{3,5\}$.)

Let $\vp$ and $\vq$ be ideals of $K$ lying above $p$ and $q$ and  denote $v$ and $w$ the corresponding places of $K$.  Since $p$ and $q$ are coprime with $\G$, the ideals $\vp$ and $\vq$ do not divide the conductor of $A$. (This is clear for if $\G = \F$. Otherwise,  $\F =1$, and, if $\vp$ or $\vq$ divide $\F_{A/K}$, then $\F \geq N_{K/\rat}(\vp)$ or $\F \geq N_{K/\rat}(\vq)$, which are both greater than 3, and this is not possible.) Hence $A$ has good reduction at $\vp$ and $\vq$ (\cite[Theorem 1]{serre-tate}). Denote $A_{v}$ and $A_{w}$ the reduced  varieties and $k_{v}$ and $k_{w}$  the residual fields. Then using the injection
$$
A(K)_{\tors} \hookrightarrow A_{v}(k_{v}) \times A_{w}(k_{w})
$$
we deduce that 
$|A(K)_{\tors}| \leq (N_{K/\rat}(\vp) \cdot N_{K/\rat}(\vq))^{g} \leq (pq)^{g[K:\rat]}
 \leq  (\frac{5}{\log 2} \log \G)^{2 g[K:\rat]}$.
We proceed in the same way for $|\Ad(K)_{\tors}|$. Since the conductor of $\Ad$ is the same as the conductor of $A$ (\cite[Corollary 2]{serre-tate}), we can conclude.
\hfill$\Box$

\subsection{Bound for the product of the order of the Tate-Shafarevic group and the regulator}\label{section-bounds-product-Sha.Reg}

We first prove Proposition \ref{prop-sha.reg}.

\medskip
\noindent {\it Proof of Proposition \ref{prop-sha.reg}.} We start by the formula (\ref{formula-bsd}) of Conjecture  \ref{bsd}. We then  bound $|L^{\star}(A/K, 1)|$ using (\ref{bound-coeff-cond})  of Lemma \ref{coeff-dominant}, the local factors using  the bound (\ref{local-falt-ineg-weaker})  of Lemma \ref{lemme-local-falt-ineg} and the torsion part of $A(K)$ using Lemma \ref{lemma-torsion}. Thus, the product  $|\Sha(A/K)| \times  \Reg(A(K)$ is bounded from above by 
\begin{equation} \label{borne-sha.reg-detailed}
(2^{16}  \cdot g^2 d)^{\frac{gd}{2}} \cdot 2^{r} \cdot D_{K}^{g} \cdot \F^{\frac{1}{4}} \cdot (\log \G)^{4gd} \cdot 
 (\log (\F \cdot  D_{K}^{2g}))^{2gd} \cdot  e^{dh} \cdot \max\{1,h\}^{\frac{dg}{2}},
\end{equation} 
which could be written as in the statement of the proposition, when $\F \ne1$. 
\hfill$\Box$\medskip

In the one-dimensional case, let's denote $E$ the abelian variety which is an elliptic curve. Then the following bounds also hold. 

\begin{prop}\label{prop-sha.reg-g=1} 

Under Hypothesis \ref{hyp-general}, the product of the order of the Tate-Shafarevic group and the regulator of the elliptic curve $E$ satisfy the following bounds
\begin{equation} \label{borne-sha.reg-ellip}
|\Sha(E/K)| \cdot \Reg(E/K) \leq C_{d} \cdot 2^r \cdot D_K \cdot \F^{\frac{1}{4}} \cdot (\log (\F \cdot D_K^2))^{2d} \cdot e^{dh} \cdot h^{d/2}, 
\end{equation}
with $C_{d} = e \left(\frac{6 \sqrt{3}}{5}\right)^{d} \cdot d^{\frac{d}{2}} \cdot (129. (5^d-1)(3d)^6)^{\frac{(1+3^{d/2})^8}{\log (1+3^{d/2})}}$, and 
\begin{equation} \label{borne-sha.reg-ellip-rank}
|\Sha(E/K)| \cdot \Reg(E/K) \leq C'_{d} \cdot D_K^{\frac{3}{2}} \cdot \F^{\frac{1}{2}} \cdot e^{dh} \cdot h^{d/2}, 
\end{equation}
with $C'_{d} = \left(\frac{9\sqrt3}{2\pi}\right)^{d} \cdot d^{\frac{d}{2}} \cdot (129. (5^d-1)(3d)^6)^{\frac{(1+3^{d/2})^8}{\log (1+3^{d/2})}}$.
\end{prop}

\noindent {\it Proof.}  
We start again from the formula (\ref{formula-bsd}) of Conjecture  \ref{bsd}. To bound the local factor, we use the bound (\ref{bound-local-falt-ineg}) of Lemma \ref{lemme-local-falt-ineg}  and to bound the torsion part of the curve, we use Lemma \ref{lemma.torsion.ellip}. We conclude for both bounds with Lemma \ref{coeff-dominant} for the leading coefficient of the $L$-series of $E/K$ at $s=1$.
\hfill$\Box$\medskip
\medskip

Our bounds of Proposition \ref{prop-sha.reg} and Proposition \ref{prop-sha.reg-g=1} extend Rémond's bounds  \cite{these.remond}, valid for an elliptic curve in the case $K= \rat$. With the same notations as in (\ref{lang'sconj-Sha.Reg}), his Proposition A.2.3, of Annex A reads as follows.
{\it Let $E$ be an elliptic curve defined over $\rat$. Suppose that $E$ verifies Conjecture  \ref{bsd}. Then, for every $\epsilon >0$, there exists a constant $C_{\epsilon}$, such that
$$| \Sha(E/\rat) \cdot \Reg (E/\rat) \leq C_{\epsilon} \F^{\frac{1}{4} + \epsilon} 2^r H(E)^{\frac{1}{12} + \epsilon};$$
and there exist an absolute constant $C > 0$ such that
$$| \Sha(E/\rat) \cdot \Reg (E/\rat) \leq C\F^{1/2 + \epsilon} H(E)^{\frac{1}{12}} (\log \max\{2, H(E)  \})^{1/2}.$$
}

In \cite{bosser-surroca-brazil}, we use estimate (\ref{borne-sha.reg-ellip-rank}), which is independent of the rank. However, estimate (\ref{borne-sha.reg-ellip}) gives a better dependence on the conductor, when the rank is neglected.

\medskip


\begin{remark}
\em{
Being inspired by the Brauer-Siegel formula\footnote{Consider the family of all number fields $K$ with degree bounded by, say, $d_0$, when the discriminant $\Delta_K$ goes to infinity. Then $\Delta_K^{1/2 - \epsilon} \ll  h_K \Reg_K \ll \Delta_K^{1/2 + \epsilon}$, where $\Reg_K$ is the regulator and $h_K$ the class number.}, one would like to also have a lower bound for the product of the order of the Tate-Shafarevich group and the canonical regulator, in terms of the height of the variety.  Pacheco and Hindry  \cite{marc-amilcar-brauer-siegel} explain why the expected lower bound seems not to be the exact translation from the Brauer-Siegel formula ($H_{Falt}(A/K)^{1-\epsilon} \ll |\Sha(A/K)| \Reg(A/K)$). Nevertheless, in \cite{autissier-hindry-pazuki}, an explicit lower bound for the regulator of an elliptic curve $E$ over a number field $K$ is given (and thus holds for the product with the order of the Tate-Shafarevich group), in terms of $[K:\rat], r=\rk(E(K)), |E(K)_{\tors}|$ and the height of the modular invariant $j_E$. Their bound depends on $j_E$ as $h(j_E)^{\frac{r-4}{3}} (\log (3h(j_E)))^{\frac{2r+2}{3}}$. Notice that $h_{Falt}(E/K) \gg \max\{h(j_E), \log N_{K/\rat} \Delta_{E/K} \}$. 

If one would like  to deduce from the BSD-conjecture  a lower bound for the product of the order of the Tate-Shafarevich group and the canonical regulator, one would be confronted with \\
- the problem of estimating from above the product $\prod_v c_v$ of the local numbers at the finite places and also with \\
- the problem of giving  a lower bound for $L^{*}(A/K,1)$. 

For the local numbers $c_v$, this could be done,  when e.g. $A$ is a jacobian variety,  under Szpiro's conjecture (see \cite[Lemma 3.5]{hindry.mordell-weil}). The question for the $L$-series also seems difficult (in the case $g=1$ and $K=\rat$ one could see the proof of Theorem 2 of \cite{goldfeld-szpiro}).}
\end{remark}

%
 
%




\section{Geometry of numbers and non-torsion points}\label{section-lower-non-torsion}

We now have a bound for the product of the order of the Tate-Shafarevic group and the regulator. Recall that the regulator is built from the canonical heights of  generators of the Mordell-Weil group, and this is what we would like to bound.  In order to manage separately these quantities, we use a classical result on geometry of numbers. 

Recall that the N\'eron-Tate height $\hat{h}_{\pp}$ on $A(K)$ 
extends to a positive definite quadratic form on  $A(K) \otimes_{\Z} \R$. We will apply Minkowski's theorem on the successive minima to the lattice $A(K)/A(K)_{\tors}$, sitting inside the euclidean space $A(K) \otimes_{\Z} \R \simeq \R^{r}$, with  inner product $<,>_{\pp}$ and which satisfies $<P,P>_{\pp} = \hat{h}_{\pp} (P)$. The symmetric convex distance-function is then $\sqrt{\hat{h}_{\pp} (.)}$,  
and  the regulator $\Reg_{\pp}(A/K)$ 
is the square of the volume of the fundamental domain for the lattice.

\begin{lemma}\label{lemma-minkowskis}
We can choose a basis $\{P_1, \ldots, P_r\}$ for the  Mordell-Weil group modulo torsion satisfying $\hat{h}_{\pp}(P_1) \leq \ldots \leq \hat{h}_{\pp}(P_r)$, and
\begin{equation}\label{minkowski}
\prod_{i=1}^{r}\hat{h}_{\pp}(P_{i}) \leq (r!)^2 r^r \Reg_{\pp}(A/K)\leq (r!)^2 \left(\frac{r}{2}\right)^{r}  \deg(\phi_{\pp})^r \Reg(A/K).
\end{equation}

\end{lemma}

\noindent {\it Proof.}  We proceed as \cite[Lemma 5.1]{remond.sous2005}. Let's denote $\lambda_1 \leq \ldots  \leq\lambda_r$  the successive minima with respect to the lattice $L=A(K)/A(K)_{\tors}$,  and $B$ the unit ball. Minkowski's theorem \cite[Theorem V, Chapter VIII, section 4.3]{cassels}  gives 
$$\lambda_1 \ldots \lambda_r  \cdot \Vol (B) \leq 2^r \Vol ((A(K) \otimes_{\Z} \R)/L) = 2^r \Reg_{\pp}(A/K)^{\frac{1}{2}}.$$
Lemma 1, page 204 of {\it loc. cit.} gives us a linear independently family of non-torsion points $Q_1, \ldots, Q_r$ such that $\lambda_i = \sqrt{\hat{h}_{\pp} (Q_i)}$.
Then, by  Lemma 8 page 135 of same reference, we know that there is a basis $P_1, \ldots, P_r$ of the torsion-free part of $A(K)$ verifying, for any $j \in \{1, \ldots, k\}$,
$$\sqrt{\hat{h}_{\pp} (P_j)} \leq \max \{\sqrt{\hat{h}_{\pp} (Q_j)}, \frac{1}{2} \sum_{l=1}^j \sqrt{\hat{h}_{\pp} (Q_l)}\} \leq \max\{1, \frac{j}{2} \}  \sqrt{\hat{h}_{\pp} (Q_j)} \leq j \lambda_j.$$
We conclude using a previous argument on the volume of the unit ball $B$ (that is, $\Vol(B)  \geq \frac{2^r}{r^{r/2}}$). 
The second inequality comes from inequality (\ref{reg}).
\hfill$\Box$\medskip
\medskip

Thus, in order to bound from below the canonical regulator of the variety, %
it suffices to give a lower bound for the  
$\hat{h}_{\pp}(P_i)'s$. %
In the same way, a lower bound for the product $\prod_{i=1}^{r-1}\hat{h}_{\pp}(P_i)$ of the $(r-1)$ first heights of the generators together with an upper bound for the canonical regulator gives us an upper bound for the greatest height $\hat{h}_{\pp}(P_r)$. 
Thus, for both applications, it will be sufficient to bound from below  the smallest height $\hat{h}_{\pp}(P_{1})$.

We are then interested in lower bounds for the height of the elements of a basis of the Mordell-Weil group modulo torsion  and more generally for points of infinite order. 
We recall that a rational point of an abelian variety has N\'eron-Tate height zero if and only if it is a torsion point. 

There are two different directions on which  these kind of lower bounds are studied. 
Let $A/K$ be an abelian variety defined over a number field. Let $K'/K$ be any finite extension of the ground field $K$ and let $P$ in $A(K')$ be a non-torsion point.
In the first case, $A/K$ is fixed and the dependence on the degree $[K':K]$ is the main interest. This is a Lehmer-type problem. In \cite{bosser-surroca-brazil} lower bounds of the first kind are used. This is because $A/K$ is fixed, while the field $K'$, which is the field of rationality of the point $P$ (which comes from a covering of $\espproj^1 \setminus \{0,1,\infty\}$), varies.  
 In the second case, the accent relies in the dependence on the variety $A/K$.  
 For the second kind of bounds, there is a conjecture of Lang \cite[page 92]{lang.dio.an}:
{\it for every elliptic curve $E/K$, there is a positive number $c_{[K:\rat]}$ depending only on the degree $[K:\rat]$, such that, for all non-torsion points $P$ in $E(K)$,}
\begin{equation}\label{lower.silv}
\hat{h}_{\pp}(P) \geq c_{[K:\rat]} \cdot \log N_{K/\rat}\Delta_{E/K}.\footnote{Hindry and Silverman \cite[Theorem 0.3]{hindry-silverman-lowerbound} proved such a lower bound for all elliptic curves, with the constant $c_{[K:\rat]} $ replaced by an explicit decreasing function of the Szpiro ratio $\sigma_{E/K} = \frac{\log N_{K/\rat} \Delta_{E/K} }{\log N_{K/\rat} \mathcal{F}_{E/K}}$.  This shows that Szpiro's conjecture implies Lang's.}
\end{equation}

Silverman \cite{silverman.duke1984} proved  Lang's conjecture for elliptic curves with integral $j$-invariant and 
 generalised it to higher dimension \cite{silverman.duke1984}: 
{\it Let $A/K$ be an abelian variety of dimension $g$, then, there exists $ c_{K,g}$ depending at most on $K$ and $g$, such that for all point $P$ in $A(K)$ generating $A$, we have
$$\hat{h}_{\pp}(P) \geq c_{K,g} h_{Falt}(A/K).
$$
}

Concerning this problem, Masser \cite[Corollary 1]{masser.LNM1290} proved that \textit{for every $K/K_{0}$, there exists a real number $c_{[K:\rat]}$ depending on $[K:\rat]$  such that, for all non-torsion points $P$ in $A(K)$, one has}
\begin{equation} \label{borne.lang.masser}
\hat{h}_{\pp}(P) \geq c_{[K:\rat]}\cdot h_{Falt}(A/K_{0})^{-(2g+1)}.
\end{equation}
We used this bound  in a first version of this work. However, Masser's result holdss for an open subset of a family of abelian varieties and the number $c_{[K:\rat]}$ is not explicit. In a second version of our work, we used Theorem 1.4 of David \cite{david.bull.soc.93}. This theorem is stated for $A_{\tau}$ the abelian variety given by the theta function $\Theta_{\tau}$ associated to $\tau$ in the Siegel space, is stated with the theta height of $A_{\tau}$, and gives a  lower bound for non-torsion points  when the abelian variety is simple. Then, using Masser-Wüstholz isogenies theorems \cite[ Lemma 4.2]{masser-wuestholz-periods1993}, we deduced the following corollary of David's result. 

\begin{prop}
Let $(A, \mathcal{L})$  be a principally polarised abelian variety of dimension $g$ defined over $K$.   Suppose that $A$ is simple  and that $A[2] \subset A(K)$.  Let $h_F(A) = \max \{ 1,h_{Stab}(A)\}$, $d =[K:\rat]$ and $D(g,d) = 3^{16g^4}\cdot 32^{4g^2}\cdot  d$. Then there exists an explicitly calculable constant $c(g)>0$ depending at most on $g$ such that for any non-torsion point $P$ in $A(K)$,
\begin{equation} \label{borne.lang.sinnou}
\hat{h}_{\mathcal{L}}(P) \geq c(g) \cdot \left(8 D(g,d)^2\right)^{-4g-2} \cdot h_F(A)^{-4g-1} \cdot  \left(\log h_F(A) + 3 \log D(g,d) \right)^{-4g-1}. 
\end{equation}
\end{prop}

In a further version of our work, using Rémond-Gaudron's isogenies theorems \cite[Th\'eor\`eme 1.2]{gaudron-remond-polaris-isog2014}, we deduced from David's theorem, a lower bound valid for $A$ not necessarily simple and not necessarily carrying a principal polarisation (as the one given by the theta function). 

Finally, we use in this present version the theorem of Bosser and Gaudron (\cite[Théorème 1.3 and Proposition 4.4]{bosser-gaudron-minorationNT}). In fact, they proved bounds valid for an abelian variety not necessarily simple, and completely explicit in both the abelian variety and the degree of the number field. (See also \cite{bruno-winckler-phd2015} for a result in this direction.) 

\begin{prop}[Bosser-Gaudron]\label{lower-bound-prop}
Let $(A, \mathcal{L})$ be a polarised abelian variety of dimension $g$ defined over a number field $K$ and let $P$ in $A(K)$ be a non-torsion point. If $g = 1$, then
\begin{equation}\label{borne.lang.bosser-gaudron-g=1}
\hat{h}_{\mathcal{L}}(P) ^{-1}  \leq 10^{40} [K:\rat]^7\max \{1, \log [K:\rat],  h_{Stab}(A)\}^{6}.
\end{equation}
If $g\geq 2$, then
\begin{equation}\label{borne.lang.bosser-gaudron-g>=2}
\hat{h}_{\mathcal{L}}(P) ^{-1}  \leq (736 g)^{8g^2}  [K:\rat] ^{4g+3} \deg(\phi_{\pp}) \max \{1, \log[K:\rat],  h_{Stab}(A)\}^{4 g + 2}.
\end{equation}
From these bounds they deduce the following bound, independent of the polarisation $\mathcal{L}$ and valid for all $g \geq 1$,
\begin{equation}\label{borne.lang.bosser-gaudron}
\hat{h}_{\mathcal{L}}(P) ^{-1}  \leq \max \{[K:\rat] + g^g,  h_{Stab}(A)\}^{6074 g}.
\end{equation}
\end{prop}

\section{On the generators of the Mordell-Weil group and the order of the Tate-Shafarevich group}\label{section-bounds}

%


In this last section, we give the proofs of Theorem \ref{thm-generateurs} and Theorem \ref{BSD+Sz-Sha}, as well as the specific results in the one-dimensional case. We then comment on these results. We prove Theorem  \ref{thm-generateurs} when proving the following result.

\medskip

\begin{prop}\label{thm-generateurs-detailed}
Suppose that $A/K$  satisfies Hypothesis \ref{hyp-general}. 
Then we can choose a system $\{P_{1}, \ldots, P_{r}\}$ of generators for the free-part of $A(K)$ verifying the following.
If $g = 1$, then
\begin{equation}\label{borne-generateurs-g=1}
\max_{1\leq i \leq r}  \hat{h}_{\pp}(P_r)  \leq C_{d} \cdot D_K \cdot (r!)^2 r^r (10^{40} d^7)^{r-1} \cdot \F^{\frac{1}{4}} (\log (\F \cdot D_K^2))^{2d} \cdot e^{dh} \cdot h^{\frac{d}{2}} \cdot \max \{1, \log d,  h\}^{6(r-1)},
\end{equation}
if $g\geq 2$, then
$$\max_{1\leq i \leq r}  \hat{h}_{\pp}(P_r)  \leq (2^{16}  g^2 d)^{\frac{gd}{2}} \cdot (r!)^2 r^r  \cdot \deg(\phi_{\pp}) ^{2r-1} \cdot (C'_{g,d})^{r-1} \cdot D_{K}^{g} \cdot \F^{\frac{1}{4}} \cdot (\log \G)^{4gd} \cdot (\log (\F D_K^{2g}))^{2gd} \cdot$$
\begin{equation}\label{borne-generateurs-g>=2}
  e^{dh} \cdot \max\{1,h\}^{\frac{dg}{2}} \cdot \max \{1, \log d,  h\}^{(4 g + 2)(r-1)},
  \end{equation}
 where  
 $\G = \max \{2, \F\}$,  $C_{d} = e \left(\frac{6 \sqrt{3}}{5}\right)^{d} \cdot d^{\frac{d}{2}} \cdot (129. (5^d-1)(3d)^6)^{\frac{(1+3^{d/2})^8}{\log (1+3^{d/2})}}$ and $C'_{g,d} = (736 g)^{8g^2} d^{4g+3}$. 
\end{prop}

\noindent {\it Proof of Theorem  \ref{thm-generateurs} and Proposition \ref{thm-generateurs-detailed}.} 
 We start from inequality (\ref{minkowski}), obtained from Minkowski's theorem on successive minima. We recall that  $\deg (\phi_{\pp}) = 1$ if $\pp$ is principal (e.g. if $g=1$). Remark that $h_{stab}(A) \leq h_{Falt}(A/K) =: h(A)$.  We then apply the conditional upper bound (\ref{borne-sha.reg-detailed}) of Proposition \ref{prop-sha.reg} for the canonical regulator to obtain the result for $g \geq 2$, and resp., the conditional upper bound (\ref{borne-sha.reg-ellip}) of Proposition \ref{prop-sha.reg-g=1} for $g=1$, and remark that $|\Sha (A/K)| \geq 1$.   We conclude applying Proposition \ref{lower-bound-prop} to $\hat{h}_{\pp}(P_{i})$, for $i = 1, ... , r-1$. Indeed, (\ref{borne.lang.bosser-gaudron-g=1}) gives the bound valid for $g=1$, (\ref{borne.lang.bosser-gaudron-g>=2}) gives the bound for $g\geq 2$, and (\ref{borne.lang.bosser-gaudron}) gives 
 $$\max_{1\leq i \leq r}  \hat{h}_{\pp}(P_r)  \leq (2^{16}  g^2 d)^{\frac{gd}{2}} \cdot r^r \cdot \deg(\phi_{\pp}) ^{2r-1} \cdot (C'_{g,d})^{r-1} \cdot D_{K}^{g} \cdot \F^{\frac{1}{4}} \cdot (\log \G)^{4gd} \cdot (\log (\F D_K^{2g}))^{2gd} \cdot$$
 $$e^{d h}\cdot \max \{d + g^g, h\}^{6074g(r-1) + \frac{dg}{2}},$$
 which could be writen as in Theorem  \ref{thm-generateurs}, when $\F \ne1$.
 \hfill$\Box$ \medskip

 The bound (\ref{borne-generateurs-g=1}), valid only for elliptic curves,  has better dependence on the conductor and on the height of the curve than (\ref{borne-generateurs}), which is more general.

\medskip

\begin{remark}\label{comparaison-lang}
\em{
For $K = \rat$,  the bound (\ref{borne-generateurs-g=1}) of Proposition \ref{thm-generateurs-detailed} becomes 
\begin{equation}
\hat{h}(P_r) \leq C_1\cdot  10^{40(r-1)} \cdot  (r!)^2 r^r \cdot  \mathcal{F}^{1/4} \cdot (\log (\mathcal{F}))^2 \cdot  e^{h} \cdot   h ^{1/2} \cdot  \max\{1, h\}^{6(r-1)},
\end{equation}
for an effective absolute constant $C_1$.  
This bound should be compared with Lang's conjecture (\ref{lang'sconj}). Lang obtained a factor $e^{r^2}$, which he could not reduce  to $e^r$, as he remarked in \cite[Note on p. 170]{lang.conj.dio}. Our bound gives a factor which grows with $r$ as $e^{(3r+1) \log r + 44 r}$. This is because we use Minkowski's theorem, instead of Hermite's, as Lang do.  Concerning the height of the variety, we have  a supplementary factor: $\max\{1, h\}^{1/2 + 6(r-1)}$. The factor $h^{1/2}$ (as well as $e^h$), comes from the local factor $c_{\infty}$. The factor $\max\{1, h\}^{6(r-1)}$ comes from the lower bound (\ref{borne.lang.bosser-gaudron-g=1}) for non-torsion points.  To bound the height of non-torsion points, in the one-dimensional case, Lang used his conjectural bound (\ref{lower.silv}) and compared the discriminant of the curve with its conductor.  As for  the dependence on the conductor, we obtain $\mathcal{F}^{1/4} \cdot (\log (\mathcal{F}))^2$. Contrary to this,  Lang suggested $\F^{\epsilon(\F)}\cdot (\log \F)^r$.
This is because, for bounding the leading coefficient of the $L$-function, he avoided the use of the functional equation, which he replaced by some hypothetical bound of his own, inspired by the Riemann hypothesis on the zeta function and  some analytic estimates. Notice that the bound (\ref{bound.ooe-top}) gives $r \ll \gamma \log \mathcal{F}$ (which  is optimal in the function field case).  This gives a bound for Lang's conjecture in terms of the conductor such as $e^{r^2} \cdot \F^{\epsilon(\F)}\cdot (\log \F)^r \ll  \mathcal{F}^{\epsilon(\mathcal{F}) + \gamma \log \log \mathcal{F}  + \gamma^2 \log \mathcal{F}}$, when our bound reads as $\mathcal{F}^{\epsilon '(\F)+ 4 \gamma \log \log \mathcal{F}+  \delta}$, with   $\epsilon '(\F) = 4 \frac{\log \log \F}{\log \F}$ and $\delta$ independent of $\F$. (Here, $g=1$ and $K = \rat$, thus $\F \geq 11$.)
}
\end{remark}
%


\medskip


We  also remark that we can bound the regulator from below  using Minkowski's inequality (\ref{minkowski}) and the lower bounds  for non-torsion points of Proposition \ref{lower-bound-prop}. 
This lower bound for $\Reg(A/K)$, together with the upper bound for the product $|\Sha(A/K)| \Reg(A/K)$ obtained in Proposition \ref{prop-sha.reg}, gives an 
estimate  for the order of $\Sha(A/K)$, which growths in the conductor and in the height as
$$|\Sha(A/K)| \ll \F^{\frac{1}{4} + \epsilon(\F)} \cdot e^{dh + \epsilon'(h)},$$ 
where $\epsilon(\F) \to 0$ when $\F \to \infty$, and $\epsilon'(h) \to 0$ when $h \to \infty$, and the implied constant in the symbol $\ll$ depends on $g, d, D_K, r,\deg(\phi_{\pp})$.

On one hand, even if this is not made explicit here, we would like to point out that there should be an inequality of the form  $ \F \ll e^{12h}$. For an elliptic curve, this is quite obvious because the Faltings' height is related to  the minimal discriminant ideal, which is divided by the conductor ($h \gg \max\{h(j_E), \log N_{K/\rat} \Delta_{E/K} \} \gg \log \F$).  In higher dimension, the implied constant in $\ll$ would depend at least on $g$ and $K$. With this inequality, we could deduce an upper bound for $|\Sha(A/K)|$, which is independent of $\F$ and growths in the height as
$$|\Sha(A/K)| \ll e^{(d+3)h + \epsilon''(h)}.$$

On the other hand, an inverse inequality between the height and the conductor would lead to an upper bound as a function in $\F$, $r$, $K$ and $g$. This inequality was predicted in the 80's by Szpiro \cite{szpiro-discri-cond}:  {\it Given $\epsilon >0$, there exists a constant $c_{\epsilon} >0$ such that for any elliptic curve $E$ defined over $\rat$ with minimal discriminant $\Delta_{E/\rat}$, we have $\log \left(|\Delta_{E/\rat} |^{\frac{1}{12}} \right) \leq \left(\frac{1}{2} + \epsilon\right) \log \F +  c_{\epsilon}$.\footnote{This inequality is optimal in the function field case.} 
} Over an arbitrary number field, we have Frey's version (see Conjecture 3.2 of \cite{hindry.mordell-weil}): {\it Given $\epsilon >0$, there exists a constant $c_{\epsilon} >0$ such that for any elliptic curve $E$ defined over a number field $K$, we have
\begin{equation}\label{szpiro-frey-conj}
 h_{Falt}(E/K) \leq \left(\frac{1}{2} + \epsilon\right) \log \F +  c_{\epsilon}.
\end{equation}
}  In higher dimension, we could expect the following (\cite{hindry.mordell-weil}).

\begin{conj}[Generalised Szpiro's conjecture]\label{gral.szpiro}
 Let $A$ be an abelian variety of dimension $g$ defined over a number field $K$. There exists  real numbers $c_{1}$ and $c_{2}$ depending at most on $g$ and $K$ such that
$$
 h_{Falt}(A/K) \leq c_{1} \log \F_{A/K} + c_{2}.
$$
\end{conj}

Looking at the function field analog and a theorem of  Deligne, Hindry suggest that we may take  $c_{1}= \left(\frac{g}{2} + \epsilon \right)$, for every $\epsilon >0$. Playing with restriction  of scalars, he adds:  $c_{2} = (g^{2} + \epsilon) \log D_{K} + c_{\epsilon, d}$, where $c_{\epsilon, d}$ depends only on $\epsilon$ and  $d=[K:\rat]$. 

We deduce Theorem \ref{BSD+Sz-Sha} from the following proposition.

\begin{prop}\label{BSD+Sz-Sha-detailed} Suppose that $A/K$ satisfies Hypothesis \ref{hyp-general}. 
Furthermore, suppose that $A/K$ satisfies  Szpiro's Conjecture (Conjecture \ref{gral.szpiro}). Then, for every $\epsilon >0$,
$$|\Sha(A/K)|  \leq  C'_{d,g,\epsilon}  \cdot D_{K}^{dg^2 + g + \epsilon}   \cdot (C_{d,g,\epsilon})^{6075 gr } \cdot (r!)^2 r^r  \cdot \F^{\frac{1}{4} + \frac{dg}{2} + \epsilon} \cdot $$
$$\cdot (\log \G)^{4gd} (\log (\F D_K^{2g}))^{2gd} (\log (c_{\epsilon, d} D_K \F))^{6075 gr + \frac{dg}{2}},$$
 where 
 $\G = \max\{2, \F\}$, $C_{d,g,\epsilon} = (d +\nolinebreak[4] g^g)(g^2 +\nolinebreak[4] \epsilon)$, and $C'_{d,g,\epsilon} = c_{\epsilon,d} (2^{16} g^2 d C_{d,g,\epsilon})^{\frac{dg}{2}}$.
 
If $g\geq2$, then the following bound also holds
$$|\Sha(A/K)|  \leq  C''_{d,g,\epsilon}  \cdot D_{K}^{dg^2 + g + \epsilon} \cdot (C^{(3)}_{d,g,\epsilon})^r   (r!)^2 r^r   \cdot\F^{\frac{1}{4} + \frac{dg}{2} + \epsilon} \cdot $$
$$\cdot (\log \G)^{4gd} (\log (\F D_K^{2g}))^{2gd} (\log (c_{\epsilon, d} D_K \F))^{(4g+2)r + \frac{dg}{2}},$$
where
 $C''_{d,g,\epsilon} = c'_{\epsilon, d}  (2^{16}g^2d (g^2 +\nolinebreak[4] \epsilon) ) ^{\frac{dg}{2}}$, and $C^{(3)}_{d,g,\epsilon}  = (736 g )^{8g^2} (d(g^2 + \epsilon))^{4g+3}$.

\end{prop}

Furthermore, when $\F \ne 1$, this second bound could be written as 
$$|\Sha(A/K)|  \leq C''_{d,g,\epsilon}  \cdot D_{K}^{dg^2 + g + \epsilon} \cdot (C^{(3)}_{d,g,\epsilon})^r   (r!)^2 r^r  \cdot \F^{\frac{1}{4} + \frac{dg}{2} + \epsilon + \gamma(\F)},$$ 
 where $\gamma(\F) = 4gd \frac{\log \log \G}{\log \G} + 2gd \frac{\log \log (\G \cdot D_K^{2g})}{\log \G} + ((4g+2)r + \frac{dg}{2}) \frac{\log \log (c_{\epsilon, d} D_K \G)}{\log \G}$ tends to 0 when $\F$ tends to infinity.

\medskip

\noindent {\it Proof.} 
With Proposition \ref{prop-sha.reg}  and (\ref{minkowski})  we obtain the following upper bound 
$$|\Sha(A/K)| \leq (2^{16} g^2 d)^{\frac{gd}{2}} \cdot D_{K}^{g} \cdot \F^{\frac{1}{4}} \cdot (\log \G)^{4gd} (\log (\F D_K^{2g}))^{2gd}  \cdot  e^{dh} \cdot \max\{1,h\}^{\frac{dg}{2}} \cdot (r!)^2 r^r  \prod_{i=1}^{r}\hat{h}_{\pp}^{-1}(P_{i}).
 $$
We now uses  the different bounds of Proposition \ref{lower-bound-prop}. 
Using (\ref{borne.lang.bosser-gaudron})  we obtain, for all $g \geq 1$, 
   \begin{equation}\label{borne-naive-sha}
 |\Sha(A/K)|  \leq 
(2^{16} g^2d) ^{\frac{dg}{2}} \cdot  (r!)^2 r^r  \cdot D_{K}^{g} \cdot \F^{\frac{1}{4}} \cdot  (\log \G)^{4gd} (\log (\F D_K^{2g}))^{2gd}  \cdot  e^{dh} \max\{d +g^g, h\}^{6075 gr + \frac{dg}{2}}.
  \end{equation}

Applying Conjecture \ref{gral.szpiro} with  $c_{1} = \left(\frac{g}{2} + \epsilon \right)$ and $c_{2}=(g^{2} + \epsilon) \log D_{K} + c_{\epsilon, d}$,  
 we obtain 
$$|\Sha(A/K)|  \leq 
c'_{\epsilon, d}  \cdot  (2^{16} g^2d)^{\frac{dg}{2} }  \cdot (r!)^2 r^r  \cdot D_{K}^{dg^2 + g + \epsilon}  \cdot \F^{\frac{1}{4} + \frac{dg}{2} + \epsilon} \cdot  (\log \G)^{4gd} (\log (\F D_K^{2g}))^{2gd}  \cdot
 $$
 $$ \cdot \max \{ d + g^g, \log (c_{\epsilon, d} \cdot D_K^{g^2 + \epsilon} \cdot \F^{g/2 + \epsilon}) \}^{6075 gr + \frac{dg}{2}},
 $$
 where $c'_{\epsilon, d}$ depends only on $\epsilon$ and  $d$.
 
 {Remark that  $A' = \Res^{K}_{\rat} A$, and, by  (\ref{A'-A}) and (\ref{N>10^{g}}),  $D_K^{2g} \F= \F_{A'/\rat} > 10^{gd}$ and this gives $\log (c_{\epsilon, d} \cdot D_K^{g^2 + \epsilon} \cdot \F^{g/2 + \epsilon} )\geq  \log (c_{\epsilon, d} \cdot 10^{\frac{g^2d}{2} + \epsilon})$  (which could be smaller than $d + g^g$).} However, this shows that $\log (c_{\epsilon, d} \cdot D_K^{g^2 + \epsilon} \cdot \F^{g/2 + \epsilon}) > 1$ and we could use the (rough) bound 
 $$ \max \{ d + g^g, \log (c_{\epsilon, d} \cdot D_K^{g^2 + \epsilon} \cdot \F^{g/2 + \epsilon}) \} \leq  (d + g^g) \cdot  \log (c_{\epsilon, d} \cdot D_K^{g^2 + \epsilon} \cdot \F^{g/2 + \epsilon})$$
$$\leq  (d + g^g) \cdot  \log (c_{\epsilon, d} \cdot D_K \cdot \F)^{g^2 + \epsilon} \leq  (d + g^g) (g^2 + \epsilon) \cdot  \log (c_{\epsilon, d} \cdot D_K \cdot \F).$$
Hence
$$|\Sha(A/K)|  \leq 
c'_{\epsilon, d}  \cdot  (2^{16} g^2d)^{\frac{dg}{2} } \cdot  [(d + g^g)(g^2 + \epsilon)]^{\frac{dg}{2}}  \cdot D_{K}^{dg^2 + g + \epsilon}   \cdot  [(d + g^g)(g^2 + \epsilon)]^{6075 gr } \cdot (r!)^2 r^r    \cdot 
 $$
 $$ \cdot \F^{\frac{1}{4} + \frac{dg}{2} + \epsilon} \cdot  (\log \G)^{4gd} (\log (\F D_K^{2g}))^{2gd}    \cdot  (\log (c_{\epsilon, d} \cdot D_K\cdot \F))^{6075 gr + \frac{dg}{2}}.$$
This proves the first item. Moreover, when $\F \ne 1$, this could be written as Theorem \ref{BSD+Sz-Sha}, where $\delta(\F) = 4gd \frac{\log \log \F}{\log \F} + 2gd \frac{\log \log (\F \cdot D_K^{2g})}{\log \F} + (6075 gr + \frac{dg}{2}) \frac{\log \log (c_{\epsilon, d} D_K \F)}{\log \F}$ tends to 0 when $\F$ tends to infinity.

Otherwise,  to prove the second item, for $g \geq 2$, using (\ref{borne.lang.bosser-gaudron-g>=2}) instead of  (\ref{borne.lang.bosser-gaudron}), we obtain, 
$$|\Sha(A/K)|  \leq 
(2^{16} g^2d) ^{\frac{dg}{2}} \cdot D_{K}^{g} \cdot (736 g )^{8g^2r} d^{(4g+3)r} \cdot (r!)^2 r^r  \cdot 
 $$
 $$e^{dh} \cdot \F^{\frac{1}{4}} \cdot  (\log \G)^{4gd} (\log (\F D_K^{2g}))^{2gd}    \cdot \max\{1, \log d, h\}^{(4g+2)r +\frac{dg}{2}}.$$
Applying Conjecture \ref{gral.szpiro}, we have $\log h \leq \log (c_{\epsilon, d} D_K^{g^2 +\epsilon} \F^{\frac{g}{2} + \epsilon})$. By the same reasoning as previously, we have $\log (c_{\epsilon, d} \cdot D_K^{g^2 + \epsilon} \cdot \F^{g/2 + \epsilon} )\geq  \log (c_{\epsilon, d} \cdot 10^{\frac{g^2d}{2} + \epsilon})$, which is larger than 1, and we could suppose larger than $\log d$ (by growing $c_{\epsilon, d}$ if needed), and then, 
$$\max\{1, \log d, h\} \leq  \log (c_{\epsilon, d} D_K^{g^2 +\epsilon} \F^{\frac{g}{2} + \epsilon}) \leq (g^2 + \epsilon) \log (c_{\epsilon, d} D_K \F ).$$
Finally
$$|\Sha(A/K)|  \leq 
c'_{\epsilon, d}  \cdot (2^{16}g^2d (g^2 + \epsilon) ) ^{\frac{dg}{2}} \cdot D_{K}^{dg^2 + g + \epsilon} \cdot (736 g )^{8g^2r} (d(g^2 + \epsilon))^{(4g+3)r} \cdot  (r!)^2 r^r   \cdot 
 $$
 $$ \F^{\frac{1}{4} + \frac{dg}{2} + \epsilon} \cdot  (\log \G)^{4gd} (\log (\F D_K^{2g}))^{2gd}    \cdot (\log (c_{\epsilon, d} D_K \F ))^{(4g+2)r +\frac{dg}{2}},$$
which achieves the proof of Proposition \ref{BSD+Sz-Sha-detailed}.
  \hfill$\Box$\medskip

In the specific case of $g=1$, we use Proposition \ref{prop-sha.reg-g=1}, and, since we are focusing on the dependence on the conductor, and, in particular, we will neglect  the dependence in the rank, we choose to use (\ref{borne-sha.reg-ellip}). 

\begin{prop}\label{BSD+Sz-Sha-g=1} 
Let $E/K$ be an elliptic curve defined over a number field $K$. Suppose that $E/K$ satisfies Hypothesis \ref{hyp-general} and 
Szpiro-Frey's Conjecture (\ref{szpiro-frey-conj}). Then, for every $\epsilon >0$,
$$|\Sha(E/K)|  \leq C_d e^{d c_{\epsilon}} \cdot D_{K}  \cdot (r!)^2 r^r  (10^{40} d^{13})^r \cdot \F^{\frac{1}{4} + \frac{d}{2} + \epsilon} \cdot (\log (\F D_K^{2}))^{2d} (\log (e^{c_{\epsilon}} \F^{\frac{1}{2} + \epsilon}))^{6r + \frac{d}{2}},$$
where one could take  
 $C_{d} = e \left(\frac{6 \sqrt{3}}{5}\right)^{d} \cdot d^{\frac{d}{2}} \cdot (129. (5^d-1)(3d)^6)^{\frac{(1+3^{d/2})^8}{\log (1+3^{d/2})}}$.
\end{prop}

When $\F \ne 1$, this bound could be written as $$|\Sha(E/K)|  \leq C_d e^{d c_{\epsilon}} \cdot D_{K}  \cdot (r!)^2 r^r  (10^{40} d^{13})^r \cdot \F^{\frac{1}{4} + \frac{d}{2} + \epsilon + \gamma'(\F)},$$
where  $\gamma'(\F) = 2d \frac{\log \log (\F D_K^{2})}{\log \F} + (6r + \frac{d}{2}) \frac{\log \log (e^{c_{\epsilon}} \F^{\frac{1}{2} + \epsilon})}{\log \F}$ tends to 0 when $\F$ tends to infinity.

\medskip
\noindent {\it Proof.} 
We start with  the bound (\ref{borne-sha.reg-ellip}) of Proposition \ref{prop-sha.reg-g=1}. Then, we apply Minkowski's theorem (\ref{minkowski}), and with  (\ref{borne.lang.bosser-gaudron-g=1}), we obtain 
 \begin{equation}\label{borne-naive-sha-g=1}
 |\Sha(E/K)|  \leq 
C_d  \cdot D_K \cdot(r!)^2 r^r  \cdot (10^{40} \cdot d^{7})^r  \cdot \F^{\frac{1}{4}} \cdot 
 (\log (\F \cdot  D_{K}^{2}))^{2d} \cdot  e^{dh} \cdot h^{\frac{d}{2}} \cdot \max\{1, \log d,  h\}^{6r},
  \end{equation}
with $C_{d} = e \left(\frac{6 \sqrt{3}}{5}\right)^{d} \cdot d^{\frac{d}{2}} \cdot (129. (5^d-1)(3d)^6)^{\frac{(1+3^{d/2})^8}{\log (1+3^{d/2})}}$.
Applyig (\ref{szpiro-frey-conj}),  we deduce 
$$ |\Sha(E/K)|  \leq 
C_d \cdot D_{K} \cdot (r!)^4  \cdot (10^{40} d^{7})^r \cdot \F^{\frac{1}{4}} \cdot 
 (\log (\F \cdot  D_{K}^{2}))^{2d} \cdot 
 $$
 $$
e^{d c_{\epsilon}} \F^{ \frac{d}{2}+\epsilon} \cdot (\log (e^{c_{\epsilon}} \F^{\frac{1}{2} + \epsilon}))^{\frac{d}{2}} \cdot \max\{1, \log d,  \log (e^{c_{\epsilon}} \F^{\frac{1}{2} + \epsilon})\}^{6r}.
 $$
We then use the following (rough) bounds: $\log d \leq d$, which is larger than 1, and since we could enlarge $c_{\epsilon}$ enough till have $\log (e^{c_{\epsilon}} \F^{\frac{1}{2} + \epsilon}) > 1$, we could bound the maximum by the product of the two elements, to obtain 
$$ \max\{1, \log d,  \log (e^{c_{\epsilon}} \F^{\frac{1}{2} + \epsilon})\} \leq d \log (e^{c_{\epsilon}} \F^{\frac{1}{2} + \epsilon}),$$
which achieves the proof of Proposition \ref{BSD+Sz-Sha-g=1}.
\hfill$\Box$\medskip

\begin{remark}\label{comparaison-goldfeld-szpiro}
\em{
In order to compare with Goldfeld-Szpiro's result, let's write the bound of Proposition \ref{BSD+Sz-Sha-g=1}  (for $K=\rat$, and thus $\F = \F_{E/\rat} \geq 11$) as 
$$
|\Sha(E/\rat)| = O (\F^{\frac{1}{4} + c + \gamma'(\F)}), \esp \textrm{with} \esp c= 1/2 + \epsilon,
$$
where the implied constant in the $O$, as well as the function $\gamma'$, depend on $r$ and $\epsilon$,  and  $\gamma'(\F)$ tends to $0$ when $\F$ tends to infinity. 
 (Putting $g=1$ and $K=\rat$ in Theorem \ref{BSD+Sz-Sha}  also gives $|\Sha(E/\rat)| = O(\F^{1/4 + c + \delta(\F)})$.) 
 This is the closest possible bound expected by   Goldfeld and Szpiro ($c > 1/2$, in \cite[page 75]{goldfeld-szpiro}). 
In Theorem 1 of {\it loc. cit.}, which we have  quoted in the introduction by (\ref{thm1-Goldfeld-Szpiro}), they  obtained $|\Sha(E/\rat)| \ll \F^{1/4 + c + \gamma(\F)}$, with $c= 3/2$ and $\gamma(\F)$  tends to $0$ when $\F$ tends to infinity. 
As we do, they use a lower bound for non-torsion points, in terms of the height of the variety, which we both then bound in terms of the conductor assuming Szpiro's conjecture. They deduce their lower bound from Lang's conjecture, together with Hindry-Silverman result on the ratio $\sigma_{E/K}$.   
Then, we both use Szpiro's conjecture a second time, for bounding the period. Indeed, the difference between the numbers $c$ is because they use the lower bound for the period:  $\Omega^{-1} \ll \Delta^{4} \ll e^{3h}$, where $\Delta$ is the minimal discriminant of the curve (see \cite[page 168]{goldfeld.banff1988}), while our Lemma \ref{lemme-local-falt-ineg} gives: $c_{\infty}^{-1} \ll h^{\frac{1}{2}} \cdot e^{h}$, which is sharper.
}
\end{remark}

\begin{remark}\label{ouverture-goldfeld-szpiro}
\em{In the same paper  Goldfeld and  Szpiro  proved \cite[Theorem 2]{goldfeld-szpiro} a sort of reciprocal statement. Precisely, they proved that if their conjectured bound (\ref{conj.g-sz.sha<cond}) for $|\Sha(E/\rat)|$ holds for {\it every} elliptic curve over $\rat$, then a weak version of Szpiro's conjecture holds ($|\Delta_{E/\rat}| \leq \F^{18+ \epsilon}$) for every elliptic curve defined over $\rat$. (The full Szpiro's conjecture $|\Delta_{E/\rat}| \leq \F^{6+ \epsilon}$ could be deduced assuming the Riemann hypothesis for the Rankin-Selberg zeta functions associated to modular forms of weight $3/2$.) The proof uses the BSD-conjecture for {\it all} elliptic curves over $\rat$, but just in the case of rank zero, which is a theorem.

It would be interesting to investigate if this result could still be obtained for a {\it fixed} elliptic curve, or how this result could be generalised to any number field or to higher dimension. 
}
\end{remark}

\noindent {\bf Aknowledgements}

 It's my pleasure to thank the colleagues who contributed to this paper (or one of the previous versions). My thanks to Carlo Gasbarri and Henri Darmon for pointing me to some references, to Jean-Benoît Bost, for encouraging me to write my first draft with more clarity, to David Masser, Sinnou David, Gaël Rémond and Pascal Autissier,  for sharing discussions on their results,  to Samuel Le Fourn, for other interesting discussions, and to Olivier Ramaré, for making explicit the constant of Lemma \ref{lemma-torsion}. I am further deeply grateful to Daniel Bertrand and Marc Hindry, for encouraging me to bringing to completion this work. Finally, I thank Gareth Jones and the EPSRC for partial support of my research.

\bibliography{mabiblio-mw-ts-acta-arith} 

\begin{thebibliography}{EMvdG19}

\bibitem[AHP18]{autissier-hindry-pazuki}
P.~Autissier, M.~Hindry, and F.~Pazuki.
\newblock Regulators of elliptic curves.
\newblock {\em arXiv:1805.03484v1}, 12 pp., 2018.

\bibitem[Aut13]{autissier-matrixlemma}
P.~Autissier.
\newblock Un lemme matriciel effectif.
\newblock {\em Math. Z.}, 273(1-2):355--361, 2013.

\bibitem[BCDT01]{breuil-conrad-diamond-taylor}
C.~Breuil, B.~Conrad, F.~Diamond, and R.~Taylor.
\newblock On the modularity of elliptic curves over {$\bold Q$}: wild 3-adic
  exercises.
\newblock {\em J. Amer. Math. Soc.}, 14(4):843--939 (electronic), 2001.

\bibitem[BFTP18]{BCGP2018}
G.~Boxer, Calegari F., Gee T., and V.~Pilloni.
\newblock Abelian surfaces over totally real fields are potentially modular.
\newblock {\em arXiv:1812.09269}, page 285, 21 December 2018.

\bibitem[BG19]{bosser-gaudron-minorationNT}
V.~Bosser and \'{E}. Gaudron.
\newblock Logarithmes des points rationnels des vari\'{e}t\'{e}s
  ab\'{e}liennes.
\newblock {\em Canad. J. Math.}, 71(2):247--298, 2019.

\bibitem[BK94]{brumer-kramer1994}
A.~Brumer and K.~Kramer.
\newblock The conductor of an abelian variety.
\newblock {\em Compositio Math.}, 92(2):227--248, 1994.

\bibitem[Blo80]{bloch1980}
S.~Bloch.
\newblock A note on height pairings, {T}amagawa numbers, and the {B}irch and
  {S}winnerton-{D}yer conjecture.
\newblock {\em Invent. Math.}, 58(1):65--76, 1980.

\bibitem[{Bos}96a]{zbMATH00908680}
J.-B. {Bost}.
\newblock {Intrinsic heights of stable varieties and abelian varieties.}
\newblock {\em {Duke Math. J.}}, 82(1):21--70, 1996.

\bibitem[{Bos}96b]{zbMATH01003129}
J.-B. {Bost}.
\newblock {P\'eriodes et isog\'enies des vari\'et\'es ab\'eliennes sur les
  corps de nombres [d'apr\`es D. Masser et G. W\"ustholz].}
\newblock In {\em {S\'eminaire Bourbaki. Volume 1994/95. Expos\'es 790-804}},
  pages 115--161, ex. Paris: Soci\'et\'e Ma\-th\'e\-ma\-tique de France, 1996.

\bibitem[BS14]{bosser-surroca-brazil}
V.~Bosser and A.~Surroca.
\newblock Elliptic logarithms, diophantine approximation and the {B}irch and
  {S}winnerton-{D}yer conjecture.
\newblock {\em Bull. Math. Soc. Brazil, New Series}, 45(1):1--23, 2014.

\bibitem[BS15]{Bhargava-Shankar2015}
M.~Bhargava and A.~Shankar.
\newblock Ternary cubic forms having bounded invariants, and the existence of a
  positive proportion of elliptic curves having rank 0.
\newblock {\em Ann. of Math. (2)}, 181(2):587--621, 2015.

\bibitem[BSD65]{birch-swinnerton-dyer}
B.~J. Birch and H.~P.~F. Swinnerton-Dyer.
\newblock Notes on elliptic curves. {II}.
\newblock {\em J. Reine Angew. Math.}, 218:79--108, 1965.

\bibitem[Cas97]{cassels}
J.~W.~S. Cassels.
\newblock {\em An introduction to the geometry of numbers}.
\newblock Classics in Mathematics. Springer-Verlag, Berlin, 1997.
\newblock Corrected reprint of the 1971 edition.

\bibitem[CS86]{cornell-silverman}
G.~Cornell and J.~H. Silverman, editors.
\newblock {\em Arithmetic geometry}.
\newblock Springer-Verlag, New York, 1986.
\newblock Papers from the conference held at the University of Connecticut,
  Storrs, Connecticut, July 30--August 10, 1984.

\bibitem[CW77]{coates-wiles1977}
J.~Coates and A.~Wiles.
\newblock On the conjecture of {B}irch and {S}winnerton-{D}yer.
\newblock {\em Invent. Math.}, 39(3):223--251, 1977.

\bibitem[Dav93]{david.bull.soc.93}
S.~David.
\newblock Minorations de hauteurs sur les vari\'et\'es ab\'eliennes.
\newblock {\em Bull. Soc. Math. France}, 121(4):509--544, 1993.

\bibitem[DP02]{sinnou-patrice.helvet}
S.~David and P.~Philippon.
\newblock Minorations des hauteurs normalis\'ees des sous-vari\'et\'es de
  vari\'et\'es abeliennes. {II}.
\newblock {\em Comment. Math. Helv.}, 77(4):639--700, 2002.

\bibitem[Dus99]{dusart1999}
P.~Dusart.
\newblock In\'{e}galit\'{e}s explicites pour {$\psi(X)$}, {$\theta(X)$},
  {$\pi(X)$} et les nombres premiers.
\newblock {\em C. R. Math. Acad. Sci. Soc. R. Can.}, 21(2):53--59, 1999.

\bibitem[Ell75]{ellison}
W.~J. Ellison.
\newblock {\em Les nombres premiers}.
\newblock Hermann, Paris, 1975.
\newblock En collaboration avec Michel Mend\`es France, Publications de
  l'Institut de Math\'ematique de l'Universit\'e de Nancago, No. IX,
  Actualit\'es Scientifiques et Industrielles, No. 1366.

\bibitem[EMvdG19]{Edixhoven-Moonen-Van-der-Geer}
B.~Edixhoven, B.~Moonen, and G.~van~der Geer.
\newblock {\em Abelian varieties}.
\newblock http://gerard.vdgeer.net/, 2019.

\bibitem[FLHS15]{freitas-bao-siksek2015}
N.~Freitas, B.~V. Le~Hung, and S.~Siksek.
\newblock Elliptic curves over real quadratic fields are modular.
\newblock {\em Invent. Math.}, 201(1):159--206, 2015.

\bibitem[Gau06]{gaudron-ENS-2006}
{\'E}.~Gaudron.
\newblock Formes lin\'eaires de logarithmes effectives sur les vari\'et\'es
  ab\'eliennes.
\newblock {\em Ann. Sci. \'Ecole Norm. Sup. (4)}, 39(5):699--773, 2006.

\bibitem[GL96]{goldfeld-lieman}
D.~Goldfeld and D.~Lieman.
\newblock Effective bounds on the size of the {T}ate-{S}hafarevich group.
\newblock {\em Math. Res. Lett.}, 3(3):309--318, 1996.

\bibitem[Gol90]{goldfeld.banff1988}
D.~Goldfeld.
\newblock Modular elliptic curves and {D}iophantine problems.
\newblock In {\em Number theory (Banff, AB, 1988)}, pages 157--175. de Gruyter,
  Berlin, 1990.

\bibitem[GR14a]{gaudron-remond-polaris-isog2014}
\'{E}. Gaudron and G.~R\'{e}mond.
\newblock Polarisations et isog\'{e}nies.
\newblock {\em Duke Math. J.}, 163(11):2057--2108, 2014.

\bibitem[GR14b]{gaudron-remond-periodes2014}
{\'E}.~Gaudron and G.~R{\'e}mond.
\newblock Th\'eor\`eme des p\'eriodes et degr\'es minimaux d'isog\'enies.
\newblock {\em Comment. Math. Helv.}, 89(2):343--403, 2014.

\bibitem[Gra01]{graftieaux2001}
P.~Graftieaux.
\newblock Formal groups and the isogeny theorem.
\newblock {\em Duke Math. J.}, (106(1):):81--121, 2001.

\bibitem[Gri18]{griffon-BS-Legendre-EC}
R.~Griffon.
\newblock Analogue of the {B}rauer-{S}iegel theorem for {L}egendre elliptic
  curves.
\newblock {\em J. Number Theory}, 193:189--212, 2018.

\bibitem[Gro82]{gross-BSD82}
B.~H. Gross.
\newblock On the conjecture of {B}irch and {S}winnerton-{D}yer for elliptic
  curves with complex multiplication.
\newblock In {\em Number theory related to Fermat's last theorem (Cambridge,
  Mass., 1981)}, volume~26 of {\em Progr. Math.}, pages 219--236. Birkh\"auser
  Boston, Mass., 1982.

\bibitem[GS95]{goldfeld-szpiro}
D.~Goldfeld and L.~Szpiro.
\newblock Bounds for the order of the {T}ate-{S}hafarevich group.
\newblock {\em Compositio Math.}, 97(1-2):71--87, 1995.
\newblock Special issue in honour of Frans Oort.

\bibitem[GZ86]{gross-zagier1986}
B.~H. Gross and D.~B. Zagier.
\newblock Heegner points and derivatives of {$L$}-series.
\newblock {\em Invent. Math.}, 84(2):225--320, 1986.

\bibitem[Hin07]{hindry.mordell-weil}
M.~Hindry.
\newblock Why is it difficult to compute the {M}ordell-{W}eil group?
\newblock {\em Proceedings of Diophantine Geometry at Centro Ennio de Giorgi,
  Pisa June 2005, Ed. Scuola Normale Superiore di Pisa}, pages 197--219, 2007.

\bibitem[HP16]{marc-amilcar-brauer-siegel}
M.~Hindry and A.~Pacheco.
\newblock An analogue of the {B}rauer-{S}iegel theorem for abelian varieties in
  positive characteristic.
\newblock {\em Mosc. Math. J.}, 16(1):45--93, 2016.

\bibitem[HS88]{hindry-silverman-lowerbound}
M.~Hindry and Joseph~H. Silverman.
\newblock The canonical height and integral points on elliptic curves.
\newblock {\em Invent. Math.}, 93(2):419--450, 1988.

\bibitem[HS00]{hindry-silverman}
M.~Hindry and Joseph~H. Silverman.
\newblock {\em Diophantine geometry}, volume 201 of {\em Graduate Texts in
  Mathematics}.
\newblock Springer-Verlag, New York, 2000.
\newblock An introduction.

\bibitem[Kol88]{kolyvagin1989}
V.~A. Kolyvagin.
\newblock Finiteness of {$E({\bf Q})$} and {$\Sha(E,{\bf Q})$} for a subclass
  of {W}eil curves.
\newblock {\em Izv. Akad. Nauk SSSR Ser. Mat.}, 52(3):522--540, 670--671, 1988.

\bibitem[KT03]{kato-trihan2003}
{K}. {K}ato and {F}. {T}rihan.
\newblock On the conjectures of {B}irch and {S}winnerton-{D}yer in
  characteristic {$p>0$}.
\newblock {\em Invent. Math.}, 153(3):537--592, 2003.

\bibitem[Lan78]{lang.dio.an}
S.~Lang.
\newblock {\em Elliptic curves: {D}iophantine analysis}, volume 231 of {\em
  Grundlehren der Mathematischen Wissenschaften [Fundamental Principles of
  Mathematical Sciences]}.
\newblock Springer-Verlag, Berlin, 1978.

\bibitem[Lan83]{lang.conj.dio}
S.~Lang.
\newblock Conjectured {D}iophantine estimates on elliptic curves.
\newblock In {\em Arithmetic and geometry, Pap. dedic. I. R. Shafarevich on the
  occasion of his sixtieth birthday. Edited by Michael Artin and John Tate,
  Vol. I}, volume~35 of {\em Progr. Math.}, pages 155--171. Birkh\"auser
  Boston, Boston, MA, 1983.

\bibitem[Lan91]{Lang-NTIII-Enciclopeadia}
S.~Lang.
\newblock {\em Number theory. {III}}, volume~60 of {\em Encyclopaedia of
  Mathematical Sciences}.
\newblock Springer-Verlag, Berlin, 1991.
\newblock Diophantine geometry.

\bibitem[LRS93]{lockhart-rosen-silverman}
P.~Lockhart, M.~Rosen, and J.~H. Silverman.
\newblock An upper bound for the conductor of an abelian variety.
\newblock {\em J. Algebraic Geom.}, 2(4):569--601, 1993.

\bibitem[Man71]{manin}
Juri~I. Manin.
\newblock Cyclotomic fields and modular curves.
\newblock {\em Uspehi Mat. Nauk}, 26(6(162)):7--71, 1971.

\bibitem[Mas87]{masser.LNM1290}
D.~W. Masser.
\newblock Small values of heights on families of abelian varieties.
\newblock In {\em Diophantine approximation and transcendence theory (Bonn,
  1985)}, volume 1290 of {\em Lecture Notes in Math.}, pages 109--148.
  Springer, Berlin, 1987.

\bibitem[Mer96]{merel}
L.~Merel.
\newblock Bornes pour la torsion des courbes elliptiques sur les corps de
  nombres.
\newblock {\em Invent. Math.}, 124(1-3):437--449, 1996.

\bibitem[Mes86]{mestre1986Compo}
J.-F. Mestre.
\newblock Formules explicites et minorations de conducteurs de vari\'et\'es
  alg\'ebriques.
\newblock {\em Compositio Math.}, 58(2):209--232, 1986.

\bibitem[Mil72]{milne}
J.~S. Milne.
\newblock On the arithmetic of abelian varieties.
\newblock {\em Invent. Math.}, 17:177--190, 1972.

\bibitem[MW93]{masser-wuestholz-periods1993}
D.~W. Masser and G.~W{\"u}stholz.
\newblock Periods and minimal abelian subvarieties.
\newblock {\em Ann. Math.}, (2, 137(2):):407 -- 458, 1993.

\bibitem[OT89]{ooe-top}
T.~Ooe and J.~Top.
\newblock On the {M}ordell-{W}eil rank of an abelian variety over a number
  field.
\newblock {\em J. Pure Appl. Algebra}, 58(3):261--265, 1989.

\bibitem[Par99]{parent.torsion1999}
P.~Parent.
\newblock Bornes effectives pour la torsion des courbes elliptiques sur les
  corps de nombres.
\newblock {\em J. Reine Angew. Math.}, 506:85--116, 1999.

\bibitem[PL08]{phragmen-lindeloef1908}
E.~Phragm{\'e}n and E.~Lindel{\"o}f.
\newblock Sur une extension d'un principe classique de l'analyse et sur
  quelques propri\'et\'es des fonctions monog\`enes dans le voisinage d'un
  point singulier.
\newblock {\em Acta Math.}, 31(1):381--406, 1908.

\bibitem[PT15]{patrikis-taylor2015}
S.~Patrikis and R.~Taylor.
\newblock Automorphy and irreducibility of some {$l$}-adic representations.
\newblock {\em Compos. Math.}, 151(2):207--229, 2015.

\bibitem[Raj97]{rajan}
C.S. Rajan.
\newblock On the size of the {T}ate-{S}hafarevich group of elliptic curves over
  function fields.
\newblock {\em Compositio Math.}, 105(1):29--41, 1997.

\bibitem[Ray85]{raynaud-asterisque127-1985}
M.~Raynaud.
\newblock Hauteurs et isog\'enies.
\newblock {\em Ast\'erisque}, (127):199--234, 1985.
\newblock Seminar on arithmetic bundles: the Mordell conjecture (Paris,
  1983/84).

\bibitem[R{\'e}m97]{these.remond}
G.~R{\'e}mond.
\newblock Sur des probl\`emes d'effectivit\'e en g\'eom\'etrie diophantienne.
\newblock {\em {T}h\`ese de doctorat, Universit\'e Paris 6},
  https://www-fourier.univ-grenoble-alpes.fr/~remond/these.ps, 1997.

\bibitem[R{\'e}m05]{remond.sous2005}
G.~R{\'e}mond.
\newblock Intersection de sous-groupes et de sous-vari\'et\'es. {I}.
\newblock {\em Math. Ann.}, 333(3):525--548, 2005.

\bibitem[R{\'e}m10]{remond-london2010}
G.~R{\'e}mond.
\newblock Nombre de points rationnels des courbes.
\newblock {\em Proc. Lond. Math. Soc. (3)}, 101:759--794, 2010.

\bibitem[Rob83]{robin1983}
G.~Robin.
\newblock Estimation de la fonction de {T}chebychef {$\theta $} sur le
  {$k$}-i\`eme nombre premier et grandes valeurs de la fonction {$\omega (n)$}
  nombre de diviseurs premiers de {$n$}.
\newblock {\em Acta Arith.}, 42(4):367--389, 1983.

\bibitem[Rub87]{rubin1987}
K.~Rubin.
\newblock Tate-{S}hafarevich groups and {$L$}-functions of elliptic curves with
  complex multiplication.
\newblock {\em Invent. Math.}, 89(3):527--559, 1987.

\bibitem[Sch03]{schoof-2003}
R.~Schoof.
\newblock Abelian varieties over cyclotomic fields with good reduction
  everywhere.
\newblock {\em Math. Ann.}, 325(3):413--448, 2003.

\bibitem[sem81]{sem-szpiro-pinceaux-ast86}
{\em {S}\'eminaire sur les {P}inceaux de {C}ourbes de {G}enre au {M}oins
  {D}eux}, volume~86 of {\em Ast\'erisque}.
\newblock Soci\'et\'e Math\'ematique de France, Paris, 1981.

\bibitem[Ser70]{serre.zeta}
J.-P. Serre.
\newblock {\em S\'eminaire {D}elange-{P}isot-{P}oitou. 11e ann\'ee: 1969/70.
  {T}h\'eorie des nombres. {F}asc. 2: {E}xpos\'e 19}.
\newblock Secr\'etariat Math\'ematique, Paris, 1970.
\newblock http://www.numdam.org.

\bibitem[Ser97]{serre}
J.-P. Serre.
\newblock {\em Lectures on the {M}ordell-{W}eil theorem}.
\newblock Aspects of Mathematics. Friedr. Vieweg \& Sohn, Braunschweig, third
  edition, 1997.
\newblock Translated from the French and edited by Martin Brown from notes by
  Michel Waldschmidt, With a foreword by Brown and Serre.

\bibitem[Shi94]{shimura.automorphic.book}
G.~Shimura.
\newblock {\em Introduction to the arithmetic theory of automorphic functions},
  volume~11 of {\em Publications of the Mathematical Society of Japan}.
\newblock Princeton University Press, Princeton, NJ, 1994.
\newblock Reprint of the 1971 original, Kano Memorial Lectures, 1.

\bibitem[Sil84]{silverman.duke1984}
J.~H. Silverman.
\newblock Lower bounds for height functions.
\newblock {\em Duke Math. J.}, 51(2):395--403, 1984.

\bibitem[ST61]{shimura-taniyama}
G.~Shimura and Y.~Taniyama.
\newblock {\em Complex multiplication of abelian varieties and its applications
  to number theory}, volume~6 of {\em Publications of the Mathematical Society
  of Japan}.
\newblock The Mathematical Society of Japan, Tokyo, 1961.

\bibitem[ST68]{serre-tate}
J.-P. Serre and J.~Tate.
\newblock Good reduction of abelian varieties.
\newblock {\em Ann. of Math. (2)}, 88:492--517, 1968.

\bibitem[Szp90]{szpiro-discri-cond}
L.~Szpiro.
\newblock Discriminant et conducteur des courbes elliptiques.
\newblock Number 183, pages 7--18. 1990.
\newblock S\'{e}minaire sur les Pinceaux de Courbes Elliptiques (Paris, 1988).

\bibitem[Tat66]{tate.bourbaki.bsd}
J.~Tate.
\newblock On the conjectures of {B}irch and {S}winnerton-{D}yer and a geometric
  analog.
\newblock In {\em S\'eminaire Bourbaki, Vol.\ 9}, pages Exp.\ No.\ 306,
  415--440 (1964--1966). Soc. Math. France, Paris, 1995, 1966.

\bibitem[{Tit}75]{titchmarsh-book-th-funct}
E.C. {Titchmarsh}.
\newblock {The theory of functions. 2nd ed.}
\newblock {London: Oxford University Press. X, 454 p. 5.00 (1975).}, 1975.

\bibitem[Wei67]{weil-math-ann1967}
A.~Weil.
\newblock \"{U}ber die {B}estimmung {D}irichletscher {R}eihen durch
  {F}unktionalgleichungen.
\newblock {\em Math. Ann.}, 168:149--156, 1967.

\bibitem[Wil95]{wiles1995}
A.~Wiles.
\newblock Modular elliptic curves and {F}ermat's last theorem.
\newblock {\em Ann. of Math. (2)}, 141(3):443--551, 1995.

\bibitem[{W}in15]{bruno-winckler-phd2015}
{B}. {W}inckler.
\newblock {I}ntersection arithm{é}tique et probl{è}me de {L}ehmer elliptique.
\newblock {\em PhD. Thesis Bordeaux}, page 120 pages, 2015.

\end{thebibliography}


\bigskip

\begin{quote}
Andrea Surroca Ortiz\\
{\tt andrea.surroca.o@gmail.com}\\
https://sites.google.com/view/andreasurroca

\end{quote}

\end{document}